\documentclass{amsart}

\usepackage{latexsym}
\usepackage{enumerate}
\usepackage{amssymb}
\usepackage{amsfonts}
\usepackage{amsmath}
\usepackage{mathrsfs}
\usepackage{amsthm}
\usepackage{geometry}
\usepackage[colorlinks,
linkcolor=red,
anchorcolor=blue,
citecolor=green
]{hyperref}
\usepackage{varioref}
\usepackage{hyperref}
\usepackage{cleveref}
\usepackage[all]{xy}
\usepackage{tikz}

\newcommand{\Ric}{\mathrm{Ric}}
\newcommand{\Euc}{\mathrm{Euc}}
\newcommand{\Vol}{\mathrm{Vol}}

\newcommand{\inj}{\mathrm{inj}}
\newcommand{\sq}{\backslash}

\newcommand{\eps}{\epsilon}
\newcommand{\Id}{\text{Id}}

\newcommand{\dist}{\mathrm{dist}}

\newcommand{\sphere}{{\mathbb{S}}}
\newcommand{\ball}{{\mathbb{B}}}

\newcommand{\tr}{\mathrm{tr}}

\newtheorem{thm}{Theorem}[section]

\newtheorem{prop}[thm]{Proposition}
\newtheorem{coro}[thm]{Corollary}
\newtheorem{lmm}[thm]{Lemma}
\newtheorem{question}[thm]{Question}

\theoremstyle{remark} 
\newtheorem*{rmk}{Remark}

\title{Examples of Ricci limit spaces with infinite bubbles}

\author{Shengxuan Zhou}
\email{zhoushx98@outlook.com}
\address{Mathematisches Institut der Universit\"at Bonn\\
Endenicher Allee 60\\
Bonn\\ D-53115\\ Germany}

\thanks{}
\date{}
\dedicatory{}

\begin{document}

\begin{abstract}
Let $k\geq 6$, $\lambda \in \mathbb{R} $, and $(X,h)$ be an $k$-dimensional smooth complete Riemannian manifold with $\Ric_h > \lambda $. In this paper, we construct, for each given $\epsilon >0$, a sequence of $(k+2)$-dimensional manifolds $(X\times\mathbb{S}^2 ,g_i ) \stackrel{GH}{\longrightarrow} (X_\epsilon ,d_\epsilon ) $ with $\Ric_{g_i} > \lambda $, such that $d_{GH} (X,X_\epsilon ) \leq \epsilon $, and $X_\epsilon $ has no open subset which is topologically a manifold. In particular, the sequence of Riemannian manifolds we construct have uniformly bounded topology, which improves the earlier examples of Hupp-Naber-Wang \cite{hnw1} and Zhou \cite{zhou1}.
\end{abstract}
\maketitle
\tableofcontents

\section{Introduction}
\label{introduction}

In the theory of Gromov-Hausdorff limits, a significant focus lies in describing the geometry and topology of measured Gromov-Hausdorff limit spaces as following:
\begin{eqnarray}
    (M^n_{i} ,g_i ,\nu_i ,p_i ) \stackrel{GH}{\longrightarrow} (X^k ,d ,\nu ,p ) ,\;\; \Ric_{g_i} \geq -\lambda ,\;\; \nu_i = \frac{1}{\Vol (B_1 (p_i)) } d\nu_{g_i },
\end{eqnarray}
where $k\in\mathbb{N}$ is the rectifiable dimension of $(X ,d ,\nu )$, and all $ M_i $ have the same dimension $n\in\mathbb{N}$. Note that Colding-Naber \cite{coldingnaber1} proved that the limit $(X ,d ,\nu )$ is $k$-rectifiable for a unique $k $.

In the non-collapsed case, which is to say $\Vol (B_1 (p_i)) \geq v >0 $, the works of Colding \cite[Theorem 0.1]{colding1} and Cheeger-Colding \cite[Theorem 5.9]{chco3} show that the dimension $k=n$, the measure $ \nu = c \mathcal{H}^n $ for some $c>0$, and the set of $n$-regular points $R_n (X) \subset X$ is a subset of full-$\nu $ measure, where $\mathcal{H}^n $ is the $n$-dimensional Hausdorff measure on $(X^n ,d )$, and a point $x\in X$ is called a $k$-regular point if and only if every tangent cone of $X$ at $x$ is isometric to $\mathbb{R}^k$. Then a Reifenberg type result of Cheeger-Colding \cite[Theorem 5.14]{chco3} shows that the limit space has a manifold structure away from a codimension $2$ set, which is the complement of an open neighborhood of $R_n (X)$. More recently, the work of Cheeger-Jiang-Naber \cite[Theorem 1.14]{jcwsjan1} was able to prove that in this case the codimension $2$ Hausdorff measure of the subset we removed can always be assumed locally finite. 

In the collapsed case, which is to say $\lim\limits_{i\to\infty} \Vol (B_1 (p_i)) = 0 $, Colding-Naber \cite{coldingnaber1} proved that the set of $k$-regular points $R_k (X) \subset X$ is a subset of full-$\nu $ measure. However, Pan-Wei's example \cite{panwei1} has shown that the Hausdorff dimension of $X\sq R_k (X) $ can be greater than the Hausdorff dimension of $ R_k (X) $. By using their almost $k$-splitting map, Cheeger-Colding \cite[Theorem 3.26]{chco5} proved that there exists a subset $U$ of $R_k (X)$ such that $\nu ( U ) >0 $ and $U$ is bi-Lipschitz to a subset of $\mathbb{R}^k$. However, in general, $U$ is not an open subset of $X$.

A natural question arises: if $x\in R_k (X)$, does there exist an open neighborhood of $x$ homeomorphic to an open subset of $\mathbb{R}^k$, similar to the non-collapsed case? When $k=1$, Honda provided an affirmative answer to this question \cite{honda1}. Recently, Bru\`e-Deng provided a proof for $k=2$ \cite{bd1}.

However, Hupp-Naber-Wang \cite{hnw1} showed that for any $4$-dimensional manifold with a Ricci lower bound, there exists a sequence of $6$-dimensional manifolds with a uniform Ricci lower bound such that the Gromov-Hausdorff limit is close to the given $4$-dimensional manifold, but the limit space has no open subset which is topologically a manifold. Moreover, Zhou \cite{zhou1} constructed $3$-dimensional Ricci limits satisfying this property.

Although the answer is negative in general, but we can ask the same question if there are some topological restriction. For example:

\begin{question}[{\cite[Question 1.1]{hnw1}}]
\label{naberopenpbm}
Let $(M^n_{i} ,g_i ,\nu_i ,p_i ) \stackrel{GH}{\longrightarrow} (X^k ,d ,\nu ,p ) $ with $\Ric_{g_i} \geq -\lambda $ and $\Vol (B_1 (p_i)) \to 0 $, where $k $ is the rectifiable dimension of $X$. Assume that $\sup_{i\in\mathbb{N}}|H_*(M_i;\mathbb{Z})|<\infty$. Does there exist an open subset of $X$ which is homeomorphic to an $k$-manifold?
\end{question}

In this paper, we construct a sequence of $(k+2)$-dimensional manifolds for any given $k$-dimensional manifolds $(k\geq 6)$ with Ricci lower bound such that the limit space has no open subset which is a topological manifold. The following theorem is our main theorem in this paper.

\begin{thm}\label{thmnonmfd}
Let $(X , h)$ be an $k$-dimensional smooth complete boundary free Riemannian manifold with $k\geq 6$ and $\Ric_h > \lambda \in\mathbb{R} $. Then for each $\epsilon >0$, there exists a metric space $(X_\epsilon ,d_\epsilon )$ such that
\begin{enumerate}[\rm (i).]
\item $d_{GH} (X,X_\epsilon ) \leq \epsilon $,  \\[-3mm]
\item $(X_\epsilon ,d_\epsilon ) $ is $k$-rectifiable, \\[-3mm]
\item There exists a sequence of $(k+2)$-dimensional boundary free complete manifolds $(X\times\mathbb{S}^2 ,g_i ) \stackrel{GH}{\longrightarrow} (X_\epsilon ,d_\epsilon )$ with $\Ric_{g_i} > \lambda $,  \\[-3mm]
\item For any nonempty open subset $U$ of $X_\epsilon $, the local $3$-homology group ${\rm{im}} \left( H_3 (U;\mathbb{Z}) \to H_3 (X_\epsilon ;\mathbb{Z}) \right) $ are infinitely generated. Consequently, every open subset of $X_\epsilon $ is not a topological manifold.
\end{enumerate}
\end{thm}

Theorem \ref{thmnonmfd} refines Hupp-Naber-Wang's example. Compared with the constructions of Hupp-Naber-Wang \cite{hnw1} and Zhou \cite{zhou1}, the sequences given in Theorem \ref{thmnonmfd} do not require any change to the topology of the underlying smooth manifold. In particular, this shows that topological conditions cannot provide a general regularity theorem about the manifold structure of collapsed Ricci limits.

However, a deep regularity theorem proved by Cheeger-Tian \cite{chti1} shows that, in $4$-dimensional case, the $L^2$ integral of curvature of Einstein metrics implies regularity in some sense. Ge-Jiang \cite{geji1} proved that Cheeger-Tian's theorem about $L^\infty$-norm of curvature operators cannot be directly generalized to higher dimensions. Nevertheless, we still hope that under certain suitable conditions, there may exist some kinds of regularity for the collapsing Ricci limits. For example, when the curvature operator $L^p$ is integrable and the global topology is fixed.

\begin{question}
\label{naberopenpbm_refine}
Fix an $n$-dimensional manifold $M$. Let $(M^n ,g_i ,\nu_i ,p_i ) \stackrel{GH}{\longrightarrow} (X^k ,d ,\nu ,p ) $, where $k\leq n-1 $ is the rectifiable dimension of $X$. Assume that:
\begin{itemize}
    \item The Ricci curvatures $\Ric_{g_i}$ satisfying the Einstein equation $\Ric_{g_i}=\Lambda g_i$ for some $\Lambda\in\{-1,0,1\}$,
    \item The diameter $\mathrm{diam}(M,g_i)\leq D$ for some $D>0$,
    \item The curvature operator $\mathfrak{R}(g_i)$ satisfying 
    $$\int_{M} \|\mathfrak{R}(g_i)\|_{g_i}^{p}d\nu_i \leq C $$ 
    for some constants $p\in(1,\infty)$ and $C>0$.
\end{itemize}

For which $p\in(1,\infty )$, does there always exist an open subset of $X$ that satisfies the above conditions and is homeomorphic to a $k$-dimensional manifold?
\end{question}

\begin{rmk}
By replacing $M$ with the unit sphere $B_1(p_i)$ in Problem \ref{naberopenpbm_refine} and removing the condition that the diameters are uniformly bounded, we can similarly formulate a local version of the topological regularity.
\end{rmk}

This paper is organized as follows. In Section \ref{sectionoutline}, we will provide an outline of our construction from topological and geometric aspects. Then we collect some preliminary results in Section \ref{sectionpreliminary}. Next, we will construct the local model in Section \ref{sectionlocalmodel}. The construction of $(X\times \mathbb{S}^2 ,g_i)$ and the proof of Theorem \ref{thmnonmfd} are contained in Section \ref{sectioninduction} and Section \ref{sectionproof}, respectively.

\vspace{0.2cm}

\textbf{Acknowledgement.}
The author would like to thank Aaron Naber for suggesting this question and providing several important ideas. The author would also be grateful to Ursula Hamenst\"adt for her interests in this topic. The author assumes full responsibility for the paper, which was created without the use of AI.

\section{Outline of the construction}
\label{sectionoutline}
In this section, we introduce the idea of our construction of Theorem \ref{thmnonmfd}. 

The basic idea is to construct a sequence of Riemannian fibrations
$$(M_i,g_i)=(X\times\sphere^2,h_i+f_i^2g_{\sphere^2})$$
with $k$-dimensional base spaces $(X,h_i)$ and $2$-dimensional fibers $\sphere^2$, where $k\geq 6$ is an integer. Moreover, as $i\to\infty$, the diameter of the fibers goes to $0$, and hence the Gromov-Hausdorff limit $X_\infty$ of $(M_i,g_i)$ is also the limit of $(X,h_i)$. This is similar to the constructions in \cite{hnw1} and \cite{zhou1}. We will then describe the construction of the sequence $(X,h_i)$ from topological and geometric ways.

\subsection{Topological part}
At first, we describe our construction in the topological sense. In this subsection, we don't consider the Ricci curvature.

Similar to \cite{zhou2}, our construction here is also based on shrinking a compact submanifold without boundary to a point. In the topological sense, our argument consists of the following four steps:
\begin{itemize}
\item Choose a smooth embedding $\iota:\bar{\ball}^3\to X$.
\item Contract the image of the boundary of $\bar{\ball}^3$, and the image $\iota(\bar{\ball}^3)$ becomes a space homeomorphic to $\sphere^3$.
\item In the step of contracting the boundary $\iota(\partial \bar{\ball}^3)$, repeat the above two steps together, and repeat this infinity times.
\item Thus, one can obtain our construction through a diagonal argument, as shown in Figure $1$ below.
\end{itemize}

%%%%%%%%%%%%%%%%%%%%%%%%%%%
% 作图区域%%%%%%%%%%%%%%%%%%%
% 颜色
\definecolor{darkgreen}{RGB}{35,175,75}
% 线条
\tikzset{
  redink/.style={thick,draw=red},
  greenink/.style={thick,draw=darkgreen},
  blueink/.style={thick,draw=blue},
  blackink/.style={thick,draw=black},
}

% 造平行四边形
\newcommand{\Panel}[2]{%
  \begin{scope}[shift={(#1,#2)}]
    \draw[thick] (0,0) -- (3.5,0) -- (3,3.25) -- (-0.5,3.25) -- cycle;
  \end{scope}
}

% 小椭圆{线类型}{x坐标}{y坐标}{x轴长}{y轴长}
\newcommand{\SplitEllipse}[5]{%
  \begin{scope}[shift={(#2,#3)}]
    \draw[#1,dashed] (#4,0)
      arc[start angle=0,end angle=-180,x radius=#4,y radius=#5];
    \draw[#1] (#4,0)
      arc[start angle=0,end angle=180,x radius=#4,y radius=#5];
  \end{scope}
}

% 气球{球顶中心x}{赤道y}{赤道x半径}{泡泡高}{颈y}{颈x}{粘合部y}{粘合部分半径}
\newcommand{\Balloon}[8]{%
  \pgfmathsetmacro{\ar}{#4*sqrt((2/3)*abs(#6-#3)/#3)}
  \pgfmathsetmacro{\al}{#4*sqrt((2/3)*abs(-#6+#3)/#3)}
  \pgfmathsetmacro{\br}{sqrt((2/3)*abs(#6-#3)*0.3*#3)}
  \pgfmathsetmacro{\bl}{sqrt((2/3)*abs(-#6+#3)*0.3*#3)}
  \pgfmathsetmacro{\tr}{sqrt((2/3)*abs(#8-#6)*0.3*#3)}

  \SplitEllipse{blackink}{#1}{#2-#5}{#6}{-0.2*#6}
  \draw[redink]
    ({#1-#8},{#2-#7})
    .. controls ({#1-#8+0.6*(#8-#6)},{#2-#7})
       and ({#1-#6},{#2-#5-\tr}) ..
    ({#1-#6},{#2-#5})
    .. controls ({#1-#6},{#2-#5+\bl})
       and ({#1-#3},{#2+\al}) ..
    ({#1-#3},#2)
    arc[start angle=180,end angle=360,x radius=#3,y radius=#4]
    .. controls ({#1+#3},{#2+\ar})
       and ({#1+#6},{#2-#5+\br}) ..
    ({#1+#6},#2-#5)
    .. controls ({#1+#6},{#2-#5-\tr})
       and ({#1+#8-0.6*(#8-#6)},{#2-#7}) ..
    ({#1+#8},#2-#7);
  \SplitEllipse{redink}{#1}{#2}{#3}{-0.2*#3}
}

\begin{tikzpicture}[scale=0.5]
  \path[use as bounding box] (0,-2) rectangle (26,17);
  \begin{scope}
    \Panel{5}{12}
    \Panel{5}{7}
    \Panel{5}{2}
    \Panel{10}{7}
    \Panel{10}{2}
    \Panel{15}{2}

    % 第一行气球{球顶中心x}{赤道y}{赤道x半径}{泡泡高}{颈y}{颈x}{粘合部y}{粘合部分半径}
    \Balloon{6.2}{13.5}{0.6}{-0.6}{0.6}{0.5}{1.2}{1}

    % 第二行气球{球顶中心x}{赤道y}{赤道x半径}{泡泡高}{颈y}{颈x}{粘合部y}{粘合部分半径}
    \Balloon{6.2}{8.5}{0.6}{-0.6}{0.6}{0.25}{1.2}{1}
    \Balloon{11.2}{8.5}{0.6}{-0.6}{0.6}{0.25}{1.2}{1}
    \Balloon{12.5}{8.8}{0.3}{-0.3}{0.3}{0.1}{0.6}{0.4}

    % 第三行气球{球顶中心x}{赤道y}{赤道x半径}{泡泡高}{颈y}{颈x}{粘合部y}{粘合部分半径}
    \Balloon{6.2}{3.5}{0.6}{-0.6}{0.6}{0.15}{1.2}{1}
    \Balloon{11.2}{3.5}{0.6}{-0.6}{0.6}{0.15}{1.2}{1}
    \Balloon{16.2}{3.5}{0.6}{-0.6}{0.6}{0.15}{1.2}{1}
    \Balloon{12.5}{3.8}{0.3}{-0.3}{0.3}{0.1}{0.6}{0.4}
    \Balloon{17.5}{3.8}{0.3}{-0.3}{0.3}{0.1}{0.6}{0.4}
    \Balloon{17}{4.7}{0.2}{-0.2}{0.2}{0.05}{0.3}{0.2}

     % 箭头
     \draw[->] (8,16) -- (20,4);
     
     % 写字
     \node at (20.2,12) {$\iota(\partial\bar{\ball}^3)\subset$};
     \node[text=red] at (23,12) {$\iota(\bar{\ball}^3)\subset$};
     \node at (24.6,12) {$X$};
     \node at (14,-0.5) {Figure $1$: Contract the embedded spheres $\iota(\partial\bar{\ball}^3)$};
  \end{scope}
\end{tikzpicture}

%%%%%%%%%%%%%%%%%%%%%%%%%%%%%
%%%%%%%%%%%%%%%%%%%%%%%%%%%%%%%%%

Note that after each contraction of an embedded $\partial\bar{\ball}^3$, the image $\iota(\bar{\ball}^3)$ becomes $\sphere^3$. Thus, in the homotopic sense, this is equivalent to the wedge sum of a $\sphere^3$ with the space. After each contraction, the wedge sum structure gives a new generator of the homology group $H_3(X_i)$. After applying this process infinitely, we have infinite generators of $H_3(X_\infty)$, which provides the construction we need.

\subsection{Geometric part}
Now we describe the metric structure in our construction.

As described in the previous subsection, the key to our construction is to contract some $\sphere^2$. In fact, we contract them by using warped products, similar to the construction in \cite{zhou2}.

Recall the handle decomposition of a $k$-sphere given by
$$\sphere^k=\bar{\ball}^3\times\sphere^{k-3}\cup_{\sphere^2\times\sphere^{k-3}}\sphere^2\times\bar{\ball}^{k-2}.$$
See also \cite[Subsection 1.4.5]{pp1}. Roughly, our construction relies on modifying the functions of the warped products to deform the part $\sphere^2\times{\ball}^{k-2}$ to become a cone $C(\sphere^2\times\sphere^{k-3})$. Note that if we take a segment $\mathfrak{I}$ connecting the origin to a point $\mathfrak{P}\in\sphere^{k-3}=\partial\sphere^{k-2}$, then by suitable smoothing, $\sphere^2\times\mathfrak{I}\cup \bar{\ball}^3\times\{\mathfrak{P}\}$ yields the required embedded $\bar{\ball}^3$. Note that in the sense of differential homeomorphism, 
$$X\cong X\# \sphere^k \# \sphere^k \# \cdots \# \sphere^k.$$

Our next step is to explain why it is possible to make the Ricci curvature of the connected sums positive by adding a new warped factor. The construction in this part is based on Perelman’s ``docking stations'' \cite{perel1} and was further developed by Burdick \cite{burdick1}, Reiser \cite{reiser1}, and others. Roughly, the construction consists of the following steps:

\begin{itemize}
\item Choose a small sphere $\sphere^{k-4}_\epsilon$ and consider the double suspension $\mathrm{S}(\mathrm{S}(\sphere^{k-4}_\epsilon))$ homeomorphic to $\sphere^{k-2}$. Note that we can always smooth the metrics if necessary.
\item Note that in the double suspension $\mathrm{S}(\mathrm{S}(\sphere^{k-4}_\epsilon))$, the singular set contains a curve $\mathfrak{S}$ given by the singular points of the first suspension $\mathrm{S}(\sphere^{k-4}_\epsilon)$. Remove a small sphere with its center inside the curve, and consider deform its boundary into a round sphere by attaching a cylinder $\mathfrak{C}$, while preserving the Ricci positive condition. Note that since we chose a small sphere in the previous step, we can ignore the restriction from the Bishop-Gromov relative volume comparison.
\item Crossing $\sphere^2$ and attach a $\sphere^{k-3}\times\bar{\ball}^3$ to the new boundary obtained in the previous step. By the handle decomposition mentioned earlier, we now obtain a space that is homeomorphic to the sphere $\sphere^k$.
\item Choose a minimal hypersurface $H$ homeomorphic to $\sphere^{k-1}$ such that the restriction of the metric on $H$ is also Ricci positive. Note that if the invariant subspace of an involution is a smooth hypersurface, then it is minimal.
\item Divide the entire $\sphere^k$ into two parts $U_1\cup U_2$ using $H$, and remove $U_2$.
\item In our construction, in fact, $H$ and a round sphere can be connected by a smooth family of metrics with positive Ricci curvature. Therefore, we can deform $H$ into a round sphere without losing the Ricci positivity of the entire space, as shown by Colding-Naber \cite{coldingnaber1}.
\item Since $U_1 \cap \mathfrak{S}\neq\varnothing$, we can take another small ball and use a similar argument. We can then attach a $\ball^{k-2}\times\sphere^2$ equipped with suitable warped product metrics.
\end{itemize}

%%%%%%%%%%%%%%%%%%%
%%%%%%%%%%%%%%%%%%%%%
% 作图区域
%%%%%%%%%%%%%%%%%%%5
% 公切线 {线种类}{圆1横坐标}{圆1纵坐标}{圆1半径}{圆2横坐标}{圆2纵坐标}{圆2半径}
\newcommand{\circletangents}[7]{%
  \begingroup
  \foreach \sign in {-1,1}{
    \draw[#1] ({#2+#4*cos((atan2(#6-#3,#5-#2))+\sign*(acos((#4-#7)/veclen(#5-#2,#6-#3))))},{#3+#4*sin((atan2(#6-#3,#5-#2))+\sign*(acos((#4-#7)/veclen(#5-#2,#6-#3))))})
      -- ({#5+#7*cos((atan2(#6-#3,#5-#2))+\sign*(acos((#4-#7)/veclen(#5-#2,#6-#3))))},{#6+#7*sin((atan2(#6-#3,#5-#2))+\sign*(acos((#4-#7)/veclen(#5-#2,#6-#3))))});
  }
  \endgroup
}

% 小球 {线种类}{横坐标}{纵坐标}{半径}
\newcommand{\balll}[4]{%
  \draw[#1] (#2,#3) circle[radius=#4];
  \begingroup
  \begin{scope}[shift={(#2,#3)}]
    \begin{scope}[shift={(0,0.1*#4)}]
      \draw[#1,dashed] (-0.99*#4,0)
        arc[start angle=180,end angle=0,x radius=0.99*#4,y radius=0.25*#4];
      \draw[#1] (0.99*#4,0)
        arc[start angle=0,end angle=-180,x radius=0.99*#4,y radius=0.25*#4];
    \end{scope}
  \end{scope}
  \endgroup
}

\begin{tikzpicture}
  \path[use as bounding box] (0,1) rectangle (14,9);
  \begin{scope}
    % 图一
    % 初始小球
    \balll{blackink}{10}{5}{0.3}
    % 一次suspension
    \draw[blackink] (10,5)
    ellipse[x radius=2.5,y radius=0.3];
    % 二次suspension-绿色为小球乘积
    % 中间弧
    \draw[blackink] (7.76,7) arc[start angle=150,end angle=210,x radius=2,y radius=4];
    % 外实线弧
    \draw[greenink] (7.9,6.7) arc[start angle=155,end angle=205,x radius=1.95,y radius=4];
    % 内实线弧
    \draw[greenink] (7.6,7.25) arc[start angle=145,end angle=215,x radius=2.1,y radius=4];
    %绿色小球
    \balll{greenink}{7.47}{4.95}{0.24}
    \balll{greenink}{7.51}{5.65}{0.24}
    \balll{greenink}{7.51}{4.25}{0.24}
    %蓝色小球
    \circletangents{blueink}{7.47}{4.95}{0.24}{6.5}{4.95}{0.17}
    \balll{blueink}{6.5}{4.95}{0.17}
    \circletangents{blueink}{5}{4.95}{0}{6}{4.95}{0.17}
    \balll{blueink}{6}{4.95}{0.17}
    \circletangents{blueink}{7.51}{4.25}{0.24}{6.6}{4}{0.17}
    \balll{blueink}{6.6}{4}{0.17}
    %切面
    \draw[redink] (6.33,4.95) -- (12.5,5);
    \draw[redink] (5,4.95) -- (6.17,4.95);
    % 写字
     \node at (6.5,1.5) {Figure $2$: Docking station};
     \node at (8,7.2) {$\mathfrak{S}$};
  \node at (11,7) {$\mathrm{S}(\mathrm{S}(\sphere^{k-4}_\epsilon))$};
  \node[text = darkgreen] at (9.8,7) {$\sphere^2\times$};
  \node[text = red] at (13,5) {$H$};
  \node at (12.5,4.5) {$U_1$};
  \node at (12.5,5.5) {$U_2$};
  \node[text = blue] at (6.7,4.55) {$\mathfrak{C}$};
  \end{scope}
\end{tikzpicture}

%%%%%%%%%%%%%%%%%%%%%%%
%%%%%%%%%%%%%%%%%%%%%%%

The above Figure 2 presents a graphical description of the main construction of the docking station. The key point is that Perelman’s construction implies that, near the singularity of the double suspension, we could modify the small ball to deform the boundary into a round sphere without losing the positivity of the Ricci curvature.

The following Figure 3 shows that after finishing the main construction of the docking station, we can slowly deform the hypersurface $H$ into a round sphere. In fact, this step can be derived directly from the Colding-Naber construction in \cite{coldingnaber1}.

%%%%%%%%%%%%%%%%%%%
%%%%%%%%%%%%%%%%%%%%%
% 作图区域
%%%%%%%%%%%%%%%%%%%5

\begin{tikzpicture}
  \path[use as bounding box] (0,1) rectangle (14,10);
  \begin{scope}
  \clip (4,4.94) -- (13,5.01) -- (13,0) -- (4,0) -- cycle;
    % 初始小球
    \balll{blackink}{10}{5}{0.3}
    % 一次suspension
    \draw[blackink] (10,5)
    ellipse[x radius=2.5,y radius=0.3];
    % 二次suspension-绿色为小球乘积
    % 中间弧
    \draw[blackink] (7.76,7) arc[start angle=150,end angle=210,x radius=2,y radius=4];
    % 外实线弧
    \draw[greenink] (7.9,6.7) arc[start angle=155,end angle=205,x radius=1.95,y radius=4];
    % 内实线弧
    \draw[greenink] (7.6,7.25) arc[start angle=145,end angle=215,x radius=2.1,y radius=4];
    %绿色小球
    \balll{greenink}{7.47}{4.95}{0.24}
    \balll{greenink}{7.51}{5.65}{0.24}
    \balll{greenink}{7.51}{4.25}{0.24}
    %蓝色小球
    \circletangents{blueink}{5}{4.95}{0}{7.47}{4.95}{0.24}
    \circletangents{blueink}{7.51}{4.25}{0.24}{6.6}{4}{0.17}
    \balll{blueink}{6.6}{4}{0.17}
  \end{scope}

  \begin{scope}
      %切面
    \draw[redink, dashed] (5,4.95) -- (12.5,5);
    \draw[redink, dashed] (8.675,5.27) -- (8.675,4.73);
    \draw[redink, dashed, rotate around={0.22:(8.675,5)}] (4.97,4.97)
        arc[start angle=180,end angle=0,x radius=3.77,y radius=0.3];
    \draw[redink, rotate around={0.22:(8.675,5)}] (12.51,4.97)
        arc[start angle=0,end angle=-180,x radius=3.77,y radius=0.3];
     %外接锥
    \draw[redink, rotate around={0.22:(8.675,9)}] (8.675,9)
        ellipse[x radius=4.5,y radius=0.4];
    \draw[redink] (4.97,4.97) -- (4.175,9);
    \draw[redink] (12.51,4.97) -- (13.175,9);
     % 写字
     \node at (6.5,2) {Figure $3$: Deform the links};
       \node[text = red] at (13,5) {$H$};
  \node at (12.5,4.5) {$U_1$};
  \node[text = blue] at (6.7,4.52) {$\mathfrak{C}$};
  \end{scope}
\end{tikzpicture}

After this part, the rest of the proof follows a similar approach to that of Hupp-Naber-Wang \cite{hnw1}, and Zhou \cite{zhou1}. Note that by modifying the additional warped factor, we can deform the cone over a small sphere into a almost Euclidean cone while preserving the positive Ricci lower bound. See also Lemma \ref{lmm_firstdoubleproduct} and Lemma \ref{huppnaberwanglemma}.

\section{Preliminary}
\label{sectionpreliminary}

\subsection{Warped products}
In this subsection, we recall the Ricci curvature of warped products.

The first case is called the double warped product.
 
Let $\varphi ,\psi $ be smooth nonnegative functions on $[0,\infty )$ such that $\varphi ,\psi $ are positive on $(0,\infty )$,
\begin{eqnarray} \label{conditionvarphi}
\varphi (0) =0 ,\; \varphi' (0)= 1 ,\; \varphi^{\mathrm{(even)}} (0) =0 ,
\end{eqnarray}
and
\begin{eqnarray} \label{conditionphirho}
\psi (0) >0 ,\; \psi^{\mathrm{(odd)}} (0) =0 .
\end{eqnarray}

Then we can define a Riemannian metric on $ \mathbb{R}^{p+1} \times \sphere^q $ by
$$ g_{\varphi ,\psi } (r) = dr^2 + \varphi (r)^2 g_{\sphere^p} +  \psi (r)^2 g_{\sphere^q} .$$
See also \cite[Proposition 1.4.7]{pp1} for more details.

Let $X_0 =\frac{\partial }{\partial r} $, $X_1 \in T_{\sphere^p} $, and $X_2 \in T_{\sphere^q} $ be unit vectors. Then the Ricci curvature of $g_{\varphi ,\psi }$ can be expressed as following. 

\begin{lmm}\label{lmmriccicurvature}
Let $\varphi ,\psi $ and $ g_{\varphi , \psi }$ be as above. Then the Ricci curvature tensor of $ g_{\varphi , \psi }$ can be determined by
\begin{eqnarray}
\Ric_{g_{\varphi , \psi } } \left(X_0 \right) & = & - \left( p \frac{\varphi''}{\varphi} + q \frac{\psi''}{\psi} \right) X_0 ,\\
\Ric_{g_{\varphi , \psi } } \left(X_1 \right) & = & \left[-  \frac{\varphi''}{\varphi} + (p-1)\frac{1-(\varphi')^2 }{\varphi^2} - q \frac{\varphi' \psi'}{\varphi \psi } \right] X_1 ,\\
\Ric_{g_{\varphi , \psi } } \left(X_2 \right) & = & \left[-  \frac{\psi''}{\psi} +(q-1) \frac{1-(\psi')^2 }{\psi^2} - p \frac{\varphi' \psi'}{\varphi \psi } \right] X_2 . \label{Remark2_Ricci}
\end{eqnarray}
\end{lmm}

\begin{proof}
This is standard, and one can conclude it by a straightforward calculation. See also \cite[Subsection 4.2.4]{pp1}.
\end{proof}

Now we consider two generalizations of the double warped product. The first one is given by the trivial sphere bundle over a Riemannian manifold.

\begin{lmm}\label{lmmriccicurvature2}
Consider the warped product $(M\times \sphere^n , g= g_M + f^2 g_{\sphere^n} )$. Then the Ricci curvature operator can be expressed as
\begin{eqnarray}
\Ric_{g} |_{TM} & = & \Ric_{g_M} -n \frac{\mathrm{Hess}^{g_M} f}{f} ,\\
\Ric_g |_{T\sphere^n} & = & \left( \frac{(n-1)(1-|\nabla^{g_M} f|_g^2)}{f^2} - \frac{\Delta^{g_M} f }{f} \right) f^2 g_{\sphere^n } .
\end{eqnarray}
\end{lmm}

\begin{proof}
Fix a point $(x,y)\in M\times\sphere^n$. Let $\{\tilde{e}_i\}_{i=1}^n$ be a vertical orthonormal basis near the slice $M\times\{y\}$ given by an orthonormal basis $\{e_i\}_{i=1}^n$ in a neighborhood of $y$ such that $\nabla^{g_{\sphere^n}}e_i=0$ at $y\in\sphere^n$, $i=1,\cdots,n$. Let $X_1,X_2\in T_M$ be local vector fields near $x$ such that $\nabla^{g_M}X_i=0$ at $x\in M$, $i=1,2$, and let $\tilde{X}_i$ be the horizontal lifting of $X_i$. By a straightforward calculation at $(x,y)$, we obtain
$$\nabla^g_{\tilde{e}_i}\tilde{e}_j=-\delta_{ij}f\nabla^{g_M}f,\quad \nabla^g_{\tilde{e}_i} \tilde{X}_1=\nabla^g_{\tilde{X}_1}\tilde{e}_i= \frac{X(f)}{f}\tilde{e}_i,\quad \nabla^g_{\tilde{X}_1}\tilde{X}_2=\widetilde{\nabla^{g_M}_{X_1}X_2} .$$
Then the lemma follows. See also \cite[Subsection 4.5.1]{pp1} and \cite[Exercise 4.7.16]{pp1}.
\end{proof}

As a corollary of Lemma \ref{lmmriccicurvature2}, one can conclude that $\Ric_{ g_M + f^2 g_{\sphere^n} } \geq \lambda $, if and only if for any $c\in (0,1]$, $\Ric_{ g_M + c^2 f^2 g_{\sphere^n} } \geq \lambda $. See also \cite[Remark 3.2]{hnw1}.

Another generalization of the double warped product is given by a bundle over an interval.

\begin{lmm}\label{lmmriccicurvature3}
Consider the warped product $([0,1)\times M, g= dr^2+h)$, where $r$ is the projection $[0,1)\times M\to [0,1)$, and $h=h(r)$ is a smooth family of Riemannian metrics on $M$. Assume that $\dim M=n$. Let $X_0=\frac{\partial}{\partial r}$, and let $X_1,X_2$ be unit vector fields tangent to $M$. Then the Ricci curvature can be expressed as
\begin{eqnarray}
\Ric_{g} \left(X_0,X_0 \right) & = & - \frac{1}{4}\left( 2\tr_hh''-|h'|_h^2 \right) ,\\
\Ric_{g} \left(X_0,X_1\right) & = & - \frac{1}{2}\left( \nabla^{h*}h'(X_1)-X_1(\tr_hh') \right),\\
\Ric_{g} \left(X_1,X_2\right) & = & \Ric_{h} \left(X_1,X_2\right)-\frac{1}{2}h''(X_1,X_2)+\frac{1}{2}h\left(h'^{\sharp}(X_1),h'^{\sharp}(X_2)\right)\\
& & -\frac{1}{4}(\tr_hh')h'(X_1,X_2), \nonumber
\end{eqnarray}
where $h'=\frac{\partial}{\partial r}h=\mathcal{L}_{X_0} h$, $h''=\frac{\partial^2}{\partial r^2}h$, $\nabla^{h*}h'(X_1)=\sum_i(\nabla_{e_i}^hh')(X_1,e_i)$, and $$h\left(h'^{\sharp}(X_1),h'^{\sharp}(X_2)\right)=\sum_ih'(X_1,e_i)h'(X_2,e_i).$$
Here, for each fixed $r$, $\{e_i\}_{i=1}^n$ is a local orthonormal basis of $(M,h(r))$.

Moreover, the second fundamental form of the boundary $\{0\}\times M$ with respect to the outward unit normal $-\frac{\partial}{\partial r}$ is given by
\begin{equation}
\label{formula_secondfundamentalform}
\mathbf{II}=-\frac{1}{2}h'(0^+)=-\frac{1}{2}\lim_{r\to 0^+} h'(r) .
\end{equation}
\end{lmm}

\begin{proof}
The lemma follows from the standard curvature equations. Note that in this warped product, the restriction of the Hessian of $r$ to $T_M$ satisfies $\mathrm{Hess}r=\frac{1}{2}h'$. See also \cite[Subsection 3.2.1]{pp1}.
\end{proof}

\subsection{Gluing lemma}

Now we recall the $C^0$ gluing lemma given by Perelman \cite[Section 4]{perel1}.

\begin{lmm}\label{lmmgluing}
Let $(M_i,g_i)$ be $n$-dimensional smooth Riemannian manifolds with compact boundary $\partial M_i$, $i=1,2$, and let $\iota:\partial M_1\to\partial M_2$ be a diffeomorphism. Assume that $\Ric_{g_i}>0$, $i=1,2$, and $\iota$ induces an isometry $(\partial M_1,g_1|_{\partial M_1})\to(\partial M_2,g_2|_{\partial M_2})$. Let $\mathbf{II}_i$ be the outer second fundamental form on $\partial M_i$, $i=1,2$. Further assume that $\iota^*\mathbf{II}_2+\mathbf{II}_1>0$ in the sense of smooth $2$-forms. Then for any open beighborhoods $U_i$ of $\partial M_i$ in $M_i$, $M_1\cup_\iota M_2$ admits a smooth Riemannian metric $\hat{g}$ such that $\Ric_{\hat{g}}>0$ and $\hat{g}|_{M_i\setminus U_i}=g_i$, $i=1,2$.
\end{lmm}

\begin{proof}
We only recall the construction, and omit the details.

Let $Z$ denote the image of $\partial M_i$ in $M=M_1\cup_\iota M_2$. Then the metric $g$ on $M$ given by $g|_{M_i}=g_i$ is Lipschitz continuous on $M$, and smooth on $M\setminus Z$. Moreover, we can identify an open neighborhood of $Z\subset U_1\cup_\iota U_2$ with $(-\epsilon,\epsilon)\times Z$, such that
\begin{itemize}
\item $Z$ correspondents to $\{0\}\times Z$,
\item $g=dr^2+h$, where $h=h(r)$ is a continuous family of smooth Riemannian metrics on $Z$,
\item the restriction of $h$ on $(-\epsilon,0]\times Z$ and $[0,\epsilon)\times Z$ are smooth tensor field.
\end{itemize}

For any $r\in(-\frac{\epsilon}{2},\frac{\epsilon}{2})$ and $\xi\in(-\frac{\epsilon}{10},\frac{\epsilon}{10})$, we can define
$$h_\xi(r)=\frac{1}{2\xi}\int_{r-\xi}^{r+\xi}h(\varrho)d\varrho.$$
Clearly, the metric $g_\xi=dr^2+h_\xi$ is $C^{1,1}$ on $(-\frac{\epsilon}{2},\frac{\epsilon}{2})\times Z$ and smooth when $r\neq \pm\xi$. Since $\iota^*\mathbf{II}_2+\mathbf{II}_1>0$, we can conclude that $h'(0^+)<h'(0^-)$. Combining this with the strictly positivity of $\Ric_g$, one can conclude that for sufficiently small $\xi$, the Ricci curvature $\Ric_{g_\xi}$ is also strictly positive. By a standard mollification of $h_\xi$ on a small open neighborhood of $r=\pm\xi$, we can get a smooth Riemannian metric $\hat{g}_\xi=dr^2+\hat{h}_\xi$ on $(-\frac{\epsilon}{2},\frac{\epsilon}{2})\times Z$ such that $\Ric_{\hat{g}_\xi}>0$ and $\hat{g}_\xi=g_\xi$ outside a small open neighborhood of $r=\pm\xi$. 

Let $\eta\in C^\infty(\mathbb{R};\mathbb{R})$ be a cut-off function such that $0\leq\eta\leq1$, $\eta=1$ on $(-\frac{\epsilon}{4},\frac{\epsilon}{4})$, and $\eta=0$ outside $(-\frac{\epsilon}{3},\frac{\epsilon}{3})$. Then for surfficiently small $\xi>0$, the metric 
$$\hat{g}=dr^2+\eta(r) \hat{h}_\xi(r)+(1-\eta(r))h(r)$$ 
is just the metric we need. This completes the proof.
\end{proof}

\subsection{Deform the links of asymptotic cones}
Note that we can use the following construction given by Colding-Naber \cite{coldingnaber1}.

\begin{lmm}[{\cite[Lemma 2.1]{coldingnaber1}}]
\label{lmmcoldingnaberexamplelemma}
Let $n\geq 3$, and let $X^{n-1}$ be a smooth compact manifold with $g(s)$, a smooth family of metrics $s\in [0,1] $ such that
\begin{enumerate}[\rm (i)${}_0$.]
\item $\Ric_{g(s)} \geq (n-2) g(s) $.
\item $\frac{d}{ds} dv(g(s)) =0 $, where $dv$ is the associated volume form.
\end{enumerate}
Then for any pair $(\lambda_1,\lambda_2)$ such that $1> \lambda_1>\lambda_2>0$, there exist a constant $R>0$ and a metric $d$ on $(C(X) ,o) = ([0,\infty )\times X /\sim ,0) $ given by $g_{\mathrm{def}}=dr^2 + f(r)^2 g(h(r)) $, such that $B_1 (o) \subset (C(X) ,d,o) $ is isometric to the unit ball $B_1 (o) \subset (C(X) ,dr^2 + \lambda_1 r^2 g(0) ,o) $, and the annulus $	A_{R,\infty } (o) \subset (C(X) ,d,o) $ is isometric to the annulus $	A_{R,\infty } (o) \subset (C(X) ,dr^2 + \lambda_2 r^2 g(1) ,o) $. 
\end{lmm}

\begin{proof}
The proof of Lemma \ref{lmmcoldingnaberexamplelemma} is exactly the same as the proof of \cite[Lemma 2.1]{coldingnaber1}.
\end{proof}

If the metric has double warped product structure, then we have the following well-known result.

\begin{lmm}
\label{lmm3dimlocalmodel}
Let $ p,q\geq 2 $. Then for any $\rho>0$, $0< \epsilon \leq \frac{1}{100} $, $0<\alpha \leq \alpha_0 (p,q,\rho,\epsilon )<1 $ and $0<\delta \leq \delta_0 (p,q,\rho,\epsilon ,\alpha) $, we can find positive constants $R(p,q,\rho,\epsilon ,\alpha ) $, $C(p,q,\rho,\epsilon ,\alpha ) $ and functions $ \varphi ,\psi $ such that 
$$\varphi(0)=1,\quad \psi (0) =\rho,\quad \min\{\varphi'(0),\psi' (0)\} \geq 0,\quad \max\{\varphi'(0),\psi' (0)\} \geq \epsilon,$$
$$\varphi |_{[R,\infty )} = (1-\epsilon ) (r+C),\quad \psi |_{[R,\infty )} = \delta (r+C)^{\alpha},\quad \Ric_{g_{\varphi ,\psi}} \geq \frac{\tau }{r^2 +1},$$
where $g_{\varphi ,\psi}=dr^2+\varphi^2(r)g_{\sphere^p}+\psi^2(r)g_{\sphere^q}$, and $\tau = \tau (p,q,\epsilon ,\alpha) >0 $ is a constant.
\end{lmm}

\begin{proof}
See \cite[Lemma 4.1]{hnw1} or \cite[Lemma 4.1]{zhou1} for details.
\end{proof}

\subsection{Perturbation near boundary}
For compact Riemannian manifolds with boundary, if the outer second fundamental form is nonnegative, then one can find a conformal perturbation near the boundary such that the outer second fundamental form becomes strictly positive.

\begin{lmm}
\label{lmm_smallpertubation_positivesecond}
Let $(M,g)$ be an $n$-dimensional compact Riemannian manifold with smooth boundary $\partial M$, and let $U$ be an open neighborhood of $\partial M$. Assume that $\Ric_g>0$, and the outer second fundamental form $\mathbf{II}_{\partial M}\geq 0$. Then there exists a conformal perturbation $\tilde{g}=e^{2u}g$ such that $u|_{M\setminus U}=u|_{\partial M}=0$, the Ricci curvature $\Ric_{\tilde{g}}>0$, and the outer second fundamental form $\widetilde{\mathbf{II}}_{\partial M}> 0$.
\end{lmm}

\begin{proof}
By a straightforward calculation,
$$\Ric_{\tilde{g}}=\Ric_{g}-(n-2)\left(\mathrm{Hess}_{g}u-du\otimes du\right) -\left(\Delta_g u+(n-2)|\nabla u|_g^2\right)g ,$$
and
$$\widetilde{\mathbf{II}}_{\partial M}=e^u \mathbf{II}_{\partial M} -e^u g(\partial_{\mathbf{n}} ,\nabla u) g,$$
where $\partial_{\mathbf{n}}$ is the outer normal vector field on $\partial M$. Now we choose a smooth function $u$ on $M$ such that $|u|+|\nabla u|_g+|\mathrm{Hess}_g u|_g$ is sufficiently small on $M$, and $g(\partial_{\mathbf{n}} ,\nabla u)<0$ on $\partial M$. By the computation above, we can conclude that $\Ric_{\tilde{g}}>0$ and $\widetilde{\mathbf{II}}_{\partial M}> 0$.
\end{proof}

\subsection{Adding conical singularities}
Now we recall the lemma Hupp-Naber-Wang used in \cite{hnw1} to add conical singularities to the given manifold.

\begin{lmm}[{\cite[Lemma 5.1]{hnw1}}]
\label{huppnaberwanglemma}
Let $( M_B ,g_B )$ be an $n$-dimensional Riemannian manifold, $n\geq 3$, $p\in M_B$, and $f$ be a positive smooth function on $B_2^{n} (p) $. Assume that $\ball_3^{n} (p) $ is pre-compact in $( M_B ,g_B )$, and the warped product space $(\ball_2^{n} (p) \times \sphere^2,g)$ with metric $g = g_B + f^2 g_{\sphere^2}$ satisfying 
\begin{align} 
\Ric_g >\lambda \, ,\inj_{g_B}(p) &\geq 2\, , \label{huppnaberwangcondition1}\\
	|\mathrm{Rm}_{g_B}|_{g_B}\, ,\, |\nabla\mathrm{Rm}_{g_B}|_{g_B}\, ,\, |\nabla^2\mathrm{Rm}_{g_B}|_{g_B}&\leq \frac{1}{2} \, , \label{huppnaberwangcondition2}\\
	|\nabla_{g_B}\ln f|_{g_B}\, ,\, |\nabla^2_{g_B}\ln f|_{g_B}\, ,\, |\nabla^3_{g_B} \ln f|_{g_B} &\leq \frac{1}{2} \, .\label{huppnaberwangcondition3}
\end{align}
Write $r := \dist_{g_B}(\cdot,p)$.
Then for all choices of parameters $0 < \eps < \eps_0(n,|\lambda|)$, $0<\alpha < \alpha_0(n,\eps)$, $0 < \hat{r} < \hat r_0(n,\alpha, \lambda,\eps)$, and $0 < \hat{\delta} < \hat{\delta}_0(n,\lambda,\alpha,\epsilon,||f||_{L^\infty},\hat{r})$, there exists $0 < \delta = \Psi(\,\hat{\delta}\|f\|_\infty\mid n, \alpha, \lambda, \eps)$ and a warped product metric $\hat g = \hat{g}_B+\hat f^{\,2} g_{\sphere^2}$ such that:
\begin{enumerate}[\rm (a).]
	\item The Ricci lower bound $\Ric_{\hat g}>\lambda - C(n)\eps$ holds for $\hat r/2 \leq r \leq 2$\, ,
	\item $\hat g = g_B+\hat{\delta}^2 f^2 g_{\sphere^2}$, $ \forall r\in [1,2] $ ,
	\item $\hat g = dr^2+(1-\eps)^2 r^2 g_{\sphere^n}+\delta^2 r^{2\alpha} g_{\sphere^2}$, $\forall r\leq \hat r$ ,
	\item The identity map $\Id:(\ball_2(p),\hat{g}_B)\to (\ball_2(p),g_B)$ is $(1+2\epsilon)$-bi-Lipschitz.
\end{enumerate}
\end{lmm}

\begin{rmk}
In this paper, $\Psi \left( \epsilon_1 ,\cdots ,\epsilon_m \big| u_1 ,\cdots ,u_n \right)$ denotes a function of $\epsilon_1 ,\cdots ,\epsilon_m $, $u_1 ,\cdots ,u_n $, such that for fixed $u_i$ we have $ \lim_{\epsilon_1 ,\cdots ,\epsilon_m \to 0} \Psi =0 $.
\end{rmk}

\begin{proof}
This lemma is a slight generalization of Hupp-Naber-Wang's lemma \cite[Lemma 5.1]{hnw1}. Although they chose $n=4$ in their original statement, the proof is essentially verbatim.  
\end{proof}

\section{Construction of the local model}
\label{sectionlocalmodel}
In this section, we construct the local model. The following proposition contains the properties we need when constructing the local model.

\begin{prop}
\label{prop_eightdim_localmodel}
Let $ k\geq 6 $. Then for any $0< \epsilon \leq \epsilon_0(k) $, $0<\alpha \leq \alpha_0 (k,\epsilon )<1 $ and $0<\delta< \delta_0 (k,\epsilon,\alpha )$, there exist:
\begin{enumerate}[\rm (i${}_{1}$).]
\item A precompact open subset $U\subset\mathbb{R}^k$,
\item A sequence of warped product metrics $\{g_j=h_j+f^2g_{\sphere^{2}}\}_{j=1}^\infty$, where $\{h_j\}_{j=1}^\infty$ is a sequence of complete Riemannian metrics on $\mathbb{R}^k$, $f$ is a smooth positive function on $\mathbb{R}^k$, and $g_{\sphere^2}$ denotes the round metric on $\sphere^2$,
\item A non-complete Riemannian manifold $(Y^\circ_{\infty},h_\infty)$,
\item A map $\mathfrak{q}:\mathbb{R}^k\to Y_\infty$, where $(Y_\infty,d_\infty)$ is the completion of $(Y^\circ_{\infty},h_\infty)$,
\item A smooth embedded ball $\iota:\bar{\ball}^3\to \mathbb{R}^k$ such that the embedded sphere $\iota(\sphere^2)=\iota(\partial \ball^3)$ contains $0^k\in\mathbb{R}^k$,
\end{enumerate}
such that
\begin{enumerate}[\rm (i${}_{1}$).]
\addtocounter{enumi}{5}
\item The remainder of $Y^\circ_{\infty}$ in its completion $Y_\infty$ is a single point $\{y_\infty\}$,
\item The sequence of pointed Riemannian manifolds $(\mathbb{R}^k,h_j,0^k)$ converges to $(Y_\infty,d_\infty,y_\infty)$ in the pointed Gromov-Hausdorff sense,
\item The map $\mathfrak{q}$ factors through the quotient $\mathbb{R}^k\to\mathbb{R}^k/\iota(\sphere^2)$, and the induced map $\mathbb{R}^k/\iota(\sphere^2)\to Y_\infty$ is a homeomorphism,
\item The restriction $\mathfrak{q}|_{\mathbb{R}^k\setminus\iota(\sphere^2)}$ gives a diffeomorphism $\mathbb{R}^k\setminus\iota(\sphere^2)\to Y^\circ_{\infty}$, and the metrics $h_j$ converging to $h_\infty$ smoothly outside $\{y_\infty\}$,
\item The restriction of $g_j$ on $(\mathbb{R}^k\setminus U)\times \sphere^2$ is independent of $j$, and it can be represented by
$$\big((\mathbb{R}^k\setminus U)\times \sphere^2,g_j|_{(\mathbb{R}^k\setminus U)\times \sphere^2}\big)\cong\left([R,\infty)\times \sphere^{k-1}\times\sphere^2,dr^2+(1-\epsilon)^2r^2g_{\sphere^{k-1}}+\delta^2r^{2\alpha}g_{\sphere^2}\right) ,$$
where $R=R(k,\epsilon,\alpha,\delta)>1$ is a positive constant,
\item The Ricci curvature $\Ric_{g_{j}}\geq \frac{\tau}{1+r_j^2}$ for some constant $\tau=\tau(k,\epsilon,\alpha,\delta)>0$, where $r_j(x)=\dist_{g_j}((0^k,*),x)$, and $*\in\sphere^2$.
\item The metrics $g_j\leq 2 g_1$, $\forall j\in\mathbb{N}$.
\end{enumerate}
\end{prop}

We divide the proof of Proposition \ref{prop_eightdim_localmodel} into three steps. First, near the singular point $y_\infty$, we construct a sequence of metrics via the warped product structure. Second, we glue this metric to a cone over a round sphere of small radius using the ``docking station" gluing construction of Perelman \cite{perel1} (see also Burdick \cite{burdick1} and Reiser \cite{reiser1}), thereby obtaining a complete metric on $\mathbb{R}^k$ with nonnegative Ricci curvature that is asymptotic to a cone. Third, applying the warped product structure again, we deform the link of this cone from a small round sphere to one close to the unit sphere, following the deformation methods of \cite{hnw1}. See also \cite{zhou1}.

\subsection{Near the singular point}
At first, we construct a sequence of complete Riemannian metrics on $\mathbb{R}^{k-2}\times\mathbb{S}^2$ with nonnegative Ricci curvature converging to a singular space. In particular, the limit space has a unique singular point of type $C(\sphere^{k-3}\times\sphere^2)$. 

\begin{lmm}
\label{lmm_firstdoubleproduct}
Let $ p,q\geq 2 $. Then there exists a sequence of complete Riemannian metrics
$$g_{p,q;i}=dr^2+\varphi_{p,q;i}^2(r)g_{\mathbb{S}^{p}}+\psi_{p,q;i}^2(r) g_{\mathbb{S}^q}$$
on $\mathbb{R}^{p+1}\times\mathbb{S}^q$, where $\varphi_{p,q;i},\psi_{p,q;i}$ are nonnegative increasing functions satisfying (\ref{conditionvarphi}) and (\ref{conditionphirho}), such that
$$ \varphi_{p,q;i}(r)=\psi_{p,q;i}(r)=e^{-10(p+q+1)}\sin(r+\epsilon_i) \quad\textrm{ for } r\in \left[\epsilon_i,\frac{\pi}{4}\right],$$
and the Ricci curvature $\Ric_{g_{p,q;i}} \geq\frac{\tau}{1+r^2} $, where $\tau=\tau(p,q)>0$ is a positive constant, and $\{\epsilon_i\}_{i=1}^\infty$ is a sequence of positive constants converging to $0$.
\end{lmm}

\begin{rmk}
By construction of $\varphi_{p,q;i}$ and $\psi_{p,q;i}$, we can conclude that for all integers $p'\in[2,p]$ and $q'\in[2,q]$, the metrics
$$g_{p',q';i}=dr^2+\varphi_{p',q';i}^2(r)g_{\mathbb{S}^{p'}}+\psi_{p',q';i}^2(r) g_{\mathbb{S}^{q'}}$$
satisfy $\Ric_{g_{p',q';i}} \geq\frac{\tau}{1+r^2} $.
\end{rmk}

\begin{proof}
Our construction is divided into two parts. At first, we construct $C^{0,1}$ functions $\hat{\varphi}_{p,q;i}$, $\hat{\psi}_{p,q;i}$. Then we will apply Lemma \ref{lmmgluing} to obtain smooth functions ${\varphi}_{p,q;i}$, ${\psi}_{p,q;i}$.

Let $\kappa_i$ denote the constant $e^{-10(p+q+i+1)}$, $\forall i\in\mathbb{Z}$. For any $i\in\mathbb{N}$, set $\xi_i=\arcsin(\frac{\sqrt{2}}{2}e^{-10i})-\frac{\kappa_i\pi}{4} $.

Define
\begin{equation}
 \hat{\varphi}_{p,q;i}(r)=\left\{
\begin{aligned}
\kappa_i\sin\left(\frac{r}{\kappa_i}\right) &,\;\;\;\; r\in \left[ 0, \frac{\pi\kappa_i}{4} \right] ,\\
 \kappa_0\sin(r+\xi_i) &,\;\;\;\; r\in \left[ \frac{\pi\kappa_i}{4},\frac{\pi}{3}-\xi_i \right] , \\
 \frac{\sqrt{3}\kappa_0}{2} \left(\frac{r}{\frac{\pi}{3}-\xi_i}\right)^{\kappa_0} &,\;\;\;\; r\in \left[ \frac{\pi}{3}-\xi_i, \infty \right) ,
\end{aligned}
\right.
\end{equation}
and
\begin{equation}
 \hat{\psi}_{p,q;i}(r)=\left\{
\begin{aligned}
\frac{\sqrt{2}}{2}\kappa_i+\kappa_0\left(\kappa_i^{-3}r^4-\frac{\pi^4\kappa_i}{256}\right) &,\;\;\;\; r\in \left[ 0, \frac{\pi\kappa_i}{4} \right] ,\\
 \kappa_0\sin(r+\xi_i) &,\;\;\;\; r\in \left[ \frac{\pi\kappa_i}{4},\frac{\pi}{3}-\xi_i \right] , \\
 \frac{\sqrt{3}\kappa_0}{2} \left(\frac{r}{\frac{\pi}{3}-\xi_i}\right)^{\kappa_0} &,\;\;\;\; r\in \left[ \frac{\pi}{3}-\xi_i, \infty \right) .\\
\end{aligned}
\right.
\end{equation}

Now we estimate the Ricci curvature for $r\notin \{0,\frac{\pi\kappa_i}{4},\frac{\pi}{3}-\xi_i\}$. Let $v_0 =\frac{\partial }{\partial r} $, $v_1 \in T_{\sphere^{k-3}} $, and $v_2 \in T_{\sphere^2} $ be unit vectors. By Lemma \ref{lmmriccicurvature}, a straightforward calculation implies that
\begin{equation}
\label{equation_Ricci_firstXzero}
 \Ric_{g_{\hat{\varphi}_{p,q;i},\hat{\psi}_{p,q;i}}}(v_0)=\left\{
\begin{aligned}
\frac{p}{\kappa_i^2}v_0-\frac{12q\kappa_0 r^2}{\frac{\sqrt{2}}{2}\kappa_i^4+\kappa_0\left(r^4-\frac{\pi^4\kappa_i^4}{256}\right)}v_0 &,\;\;\;\; r\in \left( 0, \frac{\pi\kappa_i}{4} \right) ,\\
 (p+q)v_0 &,\;\;\;\; r\in \left( \frac{\pi\kappa_i}{4},\frac{\pi}{3}-\xi_i \right) , \\
 \frac{(p+q)\kappa_0(1-\kappa_0)}{r^2}v_0 &,\;\;\;\; r\in \left( \frac{\pi}{3}-\xi_i, \infty \right) ,
\end{aligned}
\right.
\end{equation}
\begin{equation}
\label{equation_Ricci_firstXone}
 \Ric_{g_{\hat{\varphi}_{p,q;i},\hat{\psi}_{p,q;i}}}(v_1)=\left\{
\begin{aligned}
\frac{p}{\kappa_i^2}v_1-\frac{4q\kappa_0 r^3\cot\left(\frac{r}{\kappa_i}\right)}{\kappa_i\left(\frac{\sqrt{2}}{2}\kappa^4_i+\kappa_0\left(r^4-\frac{\pi^4\kappa_i^4}{256}\right)\right)}v_1 &,\;\;\;\; r\in \left( 0, \frac{\pi\kappa_i}{4} \right) ,\\
v_1 + \frac{p-1-(p+q-1)\kappa_0^2\cos^2 (r + \xi_i)}{\kappa_0^2 \sin^2 (r + \xi_i)}v_1 &,\;\;\;\; r\in \left( \frac{\pi\kappa_i}{4},\frac{\pi}{3}-\xi_i \right) , \\
 \frac{\kappa_0-(p+q)\kappa_0^2}{r^2}v_1+\frac{4(p-1)\left(\frac{\pi}{3}-\xi_i\right)^{2\kappa_0} }{3\kappa_0^2 r^{2\kappa_0}}v_1 &,\;\;\;\; r\in \left( \frac{\pi}{3}-\xi_i, \infty \right) ,
\end{aligned}
\right.
\end{equation}
and
\begin{equation}
\label{equation_Ricci_firstXtwo}
 \Ric_{g_{\hat{\varphi}_{p,q;i},\hat{\psi}_{p,q;i}}}(v_2)=\left\{
\begin{aligned}
\frac{(q-1)(\kappa_i^6 - 16 \kappa_0^2 r^6)}
{ \left( \frac{\sqrt{2}}{2} \kappa_i^4 + \kappa_0 \left(  r^4 - \frac{\pi^4 \kappa^4_i}{256} \right) \right)^{2}}v_2\quad\quad\quad\quad\quad&\\
-\frac{4 \kappa_0\kappa_i^3 r^2\left( 3\kappa_i + p r \cot \left(\frac{r}{\kappa_i}\right) \right)}{\kappa_i^4 \left( \frac{\sqrt{2}}{2} \kappa_i^4 + \kappa_0 \left(  r^4 - \frac{\pi^4 \kappa^4_i}{256} \right) \right)}
 v_2&,\;\;\;\; r\in \left( 0, \frac{\pi\kappa_i}{4} \right) ,\\
v_2 + \frac{q-1-(p+q-1)\kappa_0^2\cos^2 (r + \xi_i)}{\kappa_0^2 \sin^2 (r + \xi_i)} v_2&,\;\;\;\; r\in \left( \frac{\pi\kappa_i}{4},\frac{\pi}{3}-\xi_i \right) , \\
 \frac{\kappa_0 -(p+q)\kappa_0^2 }{r^2}v_2
+ (q-1)\frac{4\left(\frac{\pi}{3}-\xi_i\right)^{2\kappa_0} }{3\kappa_0^2 r^{2\kappa_0}} v_2&,\;\;\;\; r\in \left( \frac{\pi}{3}-\xi_i, \infty \right) .
\end{aligned}
\right.
\end{equation}

Since $\kappa_0$ is small, we can apply the formulas (\ref{equation_Ricci_firstXzero}), (\ref{equation_Ricci_firstXone}) and (\ref{equation_Ricci_firstXtwo}) to conclude that there exists a constant $\tau=\tau(p,q)>0$ satisfies 
$$\Ric_{g_{\hat{\varphi}_{p,q;i},\hat{\psi}_{p,q;i}}} \geq\frac{2\tau}{1+r^2} ,\quad \forall r\in\left( 0, \frac{\pi\kappa_i}{4} \right)\cup\left( \frac{\pi\kappa_i}{4},\frac{\pi}{3}-\xi_i \right)\cup\left( \frac{\pi}{3}-\xi_i, \infty \right). $$

By a direct computation, at $r_0\in\{0,\frac{\pi\kappa_i}{4},\frac{\pi}{3}-\xi_i\}$, we always have $\hat{\varphi}'_{p,q;i}(r_0^+)<\hat{\varphi}'_{p,q;i}(r_0^-) $ and $\hat{\psi}'_{p,q;i}(r_0^+)<\hat{\psi}'_{p,q;i}(r_0^-) $. Now we can apply Lemma \ref{lmmgluing} to construct smooth functions $\varphi_{p,q;i},\psi_{p,q;i}$ by modifying $\hat{\varphi}_{p,q;i},\hat{\psi}$ on a small open neighborhood of $\{0,\frac{\pi\kappa_i}{4},\frac{\pi}{3}-\xi_i\}$, such that the Ricci curvature $\Ric_{g_{{\varphi}_{p,q;i},{\psi}_{p,q;i}}} \geq\frac{\tau}{1+r^2} .$ Note that the second fundamental form satisfies (\ref{formula_secondfundamentalform}). Now we choose $\epsilon_i=2\kappa_i$. This completes the proof.
\end{proof}

\subsection{Docking station}
Let $k\geq 6$. In the previous step, we constructed a manifold with $\sphere^{k-2}\times\sphere^2$ as its boundary. Our next step is to introduce a suitable cobordism between $\sphere^{k-2}\times\sphere^2$ and the sphere $\sphere^{k-1}$. The following proposition provides the cobordism we need.

\begin{prop}[{\cite[Theorem 4]{burdick1}, based on \cite{perel1}}]
\label{prop_dockingstation}
For any $k\geq 6$ and $\epsilon_0\in (0,1)$, there exist a constant $R_{k,\epsilon_0}>0$ and a $k$-dimensional compact connected smooth Riemannian manifold $(\hat{Y}_{k,\epsilon_0},g_{\hat{Y}_{k,\epsilon_0}})$ with boundary $\partial{\hat{Y}}_{k,\epsilon_0}$, such that:
\begin{enumerate}[\rm (a${}_{1}$).]
\item The Ricci curvature $\Ric_{g_{\hat{Y}_{k,\epsilon_0}}}>0$,
\item The manifold $\hat{Y}_{k,\epsilon_0}$ is diffeomorphic to $\ball_2^k(0)\setminus \ball_{e^{-10k}} (\sphere_1^{2}) $, where $\ball_2^k(0)$ is the ball centered at $0\in\mathbb{R}^k$ with radius $2$, and $\sphere_1^{2}\subset\mathbb{R}^3\times\{0\}\subset \mathbb{R}^k$ is the unit sphere.
\item The boundary $\partial\hat{Y}_{k,\epsilon_0}$ has two connected components $\hat{Z}_{k,\epsilon_0;1}\sqcup \hat{Z}_{k,\epsilon_0;2}$,
\item $(\hat{Z}_{k,\epsilon_0;1},g_{\hat{Y}_{k,\epsilon_0}}|_{\hat{Z}_{k,\epsilon_0;1}})$ is isometric to $\sphere_1^{k-3}\times\sphere_1^2$ and its outer second fundamental form $\mathbf{II}_{k,\epsilon_0;1}\geq -\epsilon_0 g_{\hat{Y}_{k,\epsilon_0}}|_{\hat{Z}_{k,\epsilon_0;1}}$,
\item $(\hat{Z}_{k,\epsilon_0;2},g_{\hat{Y}_{k,\epsilon_0}}|_{\hat{Z}_{k,\epsilon_0;2}})$ is isometric to $\sphere^{k-1}_{R_{k,\epsilon_0}}$ and its outer second fundamental form $\mathbf{II}_{k,\epsilon_0;2}>0$.
\end{enumerate}
\end{prop}

The proof of Proposition \ref{prop_dockingstation}, which builds on Perelman's construction \cite{perel1}, can be found in \cite{burdick1}. See also \cite{reiser1}. For convenience, we give a simplified proof below, using ideas from Maximo-Reiser-Semola \cite{mrs1}. To prove Proposition \ref{prop_dockingstation}, 

At first, we recall two basic spaces constructed by Perelman \cite{perel1}.

\begin{lmm}
\label{lmm_Perelman_ambientandcollar}
For any $k\geq 6$, $\xi_0\in (0,e^{-10k})$ and $\epsilon\in (0,e^{-100k}\xi^{100}_0)$, there exist a Riemannian metric $g_{\mathscr{A};k,\epsilon}$ on the $(k-2)$-dimensional sphere $\sphere^{k-2}$ and a Riemannian metric $g_{\mathscr{H};k,\xi_0,\epsilon}$ on $[0,T_{\mathscr{H};k,\xi_0,\epsilon}]\times \sphere^{k-3}$, such that:
\begin{enumerate}[\rm (a${}_2$).]
\item The metric $g_{\mathscr{A};k,\epsilon}$ can be represented by
$$g_{\mathscr{A};k,\epsilon}=du^2+\cos^2 (u)d\theta^2+f^2(u) g_{\sphere^{k-4}},$$
where $u\in [0,\frac{\pi}{2}]$, $\theta\in\mathbb{R}/2\pi\mathbb{Z}$, $f|_{(0,\frac{\pi}{2}]}>0$, $f'|_{[0,\frac{\pi}{2})}>0$, $f''|_{(0,\frac{\pi}{2}]}<0$, $f'(0)=1$, $f(\frac{\pi}{2})=\epsilon$, and $f(0)=f^{(\mathrm{even})}(0)=f^{(\mathrm{odd})}(\frac{\pi}{2})=0$, 
\item On the ball $\ball_{\xi_0}(x)$ of radius $\xi_0$ centered at $x\in\{u=0\}$, the metric $g_{\mathscr{A};k,\epsilon}$ can be represented by
$$g_{\mathscr{A};k,\epsilon}|_{\ball_{\xi_0}(x)}=dr^2+\sin^2 (r)d\varsigma^2+f^2(u(r,\varsigma)) g_{\sphere^{k-4}},$$
where $\varsigma\in[-\frac{\pi}{2},\frac{\pi}{2}]$, $r(y)=\dist(x,y)$, and $u(r,\varsigma)=\arcsin(\cos(\varsigma)\sin(r))$,
\item The sectional curvature $\sec\left( g_{\mathscr{A};k,\epsilon}|_{\partial\ball_{\xi_0}(x)}\right)>\cot^2 (\xi_0)$, and the outer second fundamental form $\mathbf{II}_{\partial\ball_{\xi_0}(x)} \leq \cot (\xi_0)g_{\mathscr{A};k,\epsilon}|_{\partial\ball_{\xi_0}(x)}$,
\item The metric $g_{\mathscr{H};k,\xi_0,\epsilon}$ can be represented by
$$g_{\mathscr{H};k,\xi_0,\epsilon}=dt^2+ A^2(t,\varkappa)d\varkappa^2+ B^2(t,\varkappa)g_{\sphere^{k-4}},$$
where $t\in[0,T_{\mathscr{H};k,\xi_0,\epsilon}]$, $\varkappa\in[0,\pi]$, $A>0$, $B|_{\{\varkappa\in(0,\pi)\}}>0$,
$$\left(\sphere^{k-3},g_{\mathscr{H};k,\xi_0,\epsilon}\big|_{\{t=0\}}\right)\textrm{ is isometric to } \left(\sphere^{k-3},g_{\sphere^{k-3}}\right) ,$$
$$\left(\sphere^{k-3},g_{\mathscr{H};k,\xi_0,\epsilon}\big|_{\{t=T_{\mathscr{H};k,\xi_0,\epsilon}\}}\right)\textrm{ is isometric to } \left(\partial\ball_{\xi_0}(x),R_{\mathscr{H};k,\xi_0,\epsilon}^2g_{\mathscr{A};k,\epsilon}\big|_{\partial\ball_{\xi_0}(x)}\right) ,$$
where $R_{\mathscr{H};k,\xi_0,\epsilon}>1$ is a constant, and
$$\frac{\partial g_{\mathscr{H};k,\xi_0,\epsilon}}{\partial t}\bigg|_{t=0}\leq 2\xi_0g_{\mathscr{H};k,\xi_0,\epsilon},\quad \frac{\partial g_{\mathscr{H};k,\xi_0,\epsilon}}{\partial t}\bigg|_{t=T_{\mathscr{H};k,\xi_0,\epsilon}}> \frac{\cot(\xi_0)}{R_{\mathscr{H};k,\xi_0,\epsilon}}g_{\mathscr{H};k,\xi_0,\epsilon},$$
\item The Ricci curvature $\Ric_{g_{\mathscr{A};k,\epsilon}}>0$, and $\Ric_{g_{\mathscr{H};k,\xi_0,\epsilon}}>0$.
\end{enumerate}
\end{lmm}

\begin{proof}
At first, we need to construct suitable functions $f$ satisfying (a${}_2$).

Let $\hat{f}=\hat{f}(u)$ be a Lipschitz continuous nonnegative function on $[0,\frac{\pi}{2}]$, such that $\hat{f}$ is smooth on $[0, \frac{\pi\epsilon^{100}}{4})\cup( \frac{\pi\epsilon^{100}}{4},\epsilon^4)\cup(\epsilon^4,\frac{\pi}{6}-\zeta_{k})\cup(\frac{\pi}{6}-\zeta_{k},\frac{\pi}{2}]$, and
\begin{equation}
\left\{
\begin{aligned}
\hat{f}(u)=\epsilon^{100}\sin\left(\epsilon^{-100}u\right) &,\;\;\;\; u\in \left[ 0, \frac{\pi\epsilon^{100}}{4} \right] ,\\
\hat{f}''(u)\leq -2\epsilon^{-1}\hat{f}(u)  &,\;\;\;\; u\in \left[ \frac{\pi\epsilon^{100}}{4},\epsilon^{4} \right] ,\\
\epsilon\sin\left(u+\zeta_{k}\right)  &,\;\;\;\; u\in \left[ \epsilon^{4},\frac{\pi}{6}-\zeta_{k} \right] ,\\
\frac{\epsilon}{2}\sin\left(\frac{6u-\pi+6\zeta_{k}}{2\pi+6\zeta_{k}}\right)+\frac{\epsilon}{2}  &,\;\;\;\; u\in \left[ \frac{\pi}{6}-\zeta_{k},\frac{\pi}{2} \right] ,\\
\hat{f}(u_-)>\hat{f}'(u_+) &,\;\;\;\; u\in\left\{\frac{\pi\epsilon^{10}}{4},\frac{\pi}{6}-\zeta_{k},\epsilon^2\right\},
\end{aligned}
\right.
\end{equation}
where $\zeta_{k}\in(0,\epsilon^{20})$ is a constant. By modifying $\hat{f}$ on a small open neighborhood of $\left\{\frac{\pi\epsilon^{100}}{4},\epsilon^4,\frac{\pi}{6}-\zeta_{k}\right\}$, we can construct a function $f$ satisfying the conditions in (a${}_2$). By construction, $g_{\mathscr{A};k,\epsilon}$ is a smooth Riemannian metric, which proves (a${}_2$). Since $g_{\mathscr{A};k,\epsilon}$ is by definition a warped product over $\sphere^2$, we can conclude (b${}_2$) by a straightforward calculation.

Now we consider (c${}_2$). By construction, we have
\begin{equation}
\label{conditions_Perelmanneck_functionF}
\left\{
\begin{aligned}
f''(u)\leq \hat{f}''(u)&,\;\;\;\; u\notin\left\{\frac{\pi\epsilon^{100}}{4},\epsilon^4\right\} ,\\
f'(u)\leq \hat{f}'(u)&,\;\;\;\; u\notin\left\{\frac{\pi\epsilon^{100}}{4},\epsilon^4\right\} ,\\
f''(u)\leq -\epsilon^{-1}f(u)  &,\;\;\;\; u\in \left[ 0,\epsilon^{4}+\epsilon^{100} \right] ,\\
|f(u)-\hat{f}(u)|+|f'(u)-\hat{f}'(u)|\leq \epsilon^{10}  &,\;\;\;\; u\in \left[ \epsilon^{4}+\epsilon^{100},\frac{\pi}{8} \right] ,\\
f''(u)\leq - f(u)  &,\;\;\;\; u\in \left[ 0,\frac{\pi}{8} \right] ,\\
f'(u)\leq f(u)\cot(u)  &,\;\;\;\; u\in \left[ 0,\frac{\pi}{8} \right] .\\
\end{aligned}
\right.
\end{equation}
By a straightforward calculation,
\begin{eqnarray*}
\mathbf{II}_{\partial\ball_{\xi_0}(x)} & = & \frac{1}{2}\frac{\partial g_{\mathscr{A};k,\epsilon}}{\partial r}\bigg|_{r=\xi_0} =\sin(\xi_0)\cos(\xi_0) d\varsigma^2+f(u)f'(u)u_r(\xi_0,\varsigma) g_{\sphere^{k-4}}\\
&\leq& \max\left\{\cot(\xi_0),\frac{f'(u)}{f(u)}\cdot\frac{ \cos (\varsigma)\cos(\xi_0)}{\cos(u)}\right\}  g_{\mathscr{A};k,\epsilon}|_{\partial\ball_{\xi_0}(x)}\\
&\leq& \cot(\xi_0)   g_{\mathscr{A};k,\epsilon}|_{\partial\ball_{\xi_0}(x)} ,
\end{eqnarray*}
where $u_r(r,\varsigma)=\frac{\partial u}{\partial r}(r,\varsigma) = \frac{ \cos (\varsigma)\cos(r)}{\cos(u)} $. 

By the warped product structure on $g_{\mathscr{A};k,\epsilon}|_{\partial\ball_{\xi_0}(x)}$, the sectional curvature $\geq (1+\epsilon^{100})\cot^2(\xi_0)$ if and only if at any point $(\xi_0,\varsigma)$,
\begin{equation}
\left\{
\begin{aligned}
\frac{1-f'^2(u)\cot^2(\xi_0)\left(\tan^2(\xi_0)-\tan^2(u)\right)}{f^2(u)} \geq (1+\epsilon^{100})\cot^2(\xi_0) &, \\
\frac{-f''(u)\left(\tan^2(\xi_0)-\tan^2(u)\right)+f'(u)\left(\tan^2(u)+1\right)\tan(u)}{f(u)\tan^2(\xi_0)} \geq (1+\epsilon^{100})\cot^2(\xi_0) & .\\
\end{aligned}
\right.
\end{equation}
The two estimates also follow from (\ref{conditions_Perelmanneck_functionF}). This proves (c${}_2$).

Our next step is to construct functions $A$ and $B$ satisfying (d${}_2$).

For any $\varsigma\in[-\frac{\pi}{2},\frac{\pi}{2}]$, let $\varkappa=\varkappa(\varsigma)\in[-\frac{\pi}{2},\frac{\pi}{2}]$ defined by $\mathrm{sgn}(\varkappa)=\mathrm{sgn}(\varsigma)$ and $\cos(\varkappa)f(\xi_0)=f(u(\xi_0,\varsigma))$. It is easy to see that $\varkappa$ is smooth on $(-\frac{\pi}{2},\frac{\pi}{2})$. By a reparameterization, we have
$$\cot^2(\xi_0)g_{\mathscr{A};k,\epsilon}|_{\partial \ball_{\xi_0}(x)}=A^2_1(\varkappa)d\varkappa^2+r_1^2\cos^2(\varkappa) g_{\sphere^{k-4}}, $$
where $r_1=f(\xi_0)\cot(\xi_0)$, and $A_1>0$ be a smooth function on $\varsigma\in[-\frac{\pi}{2},\frac{\pi}{2}]$. Since the integral of $A_1$ is $\pi\cos(\xi_0)$, $\sup A_1>\cos(\xi_0)$. Moreover, by the lower bound of sectional curvature proved in (c${_2}$), we can conclude that $\sup A_1\leq 1-\epsilon^{101}$. 

Fix $\delta=e^{-\epsilon^{-100}}$ and $m=e^{\delta^{-100}}$. Let $a_\infty=\frac{\sup A_1}{r_1}$, $\alpha=(1+\delta)\log (a_\infty)$, $\beta=(1-\delta)(\log(\xi_0^2)-\log (r_1))$, $\eta(\varkappa)=\frac{A_1(\varkappa)-r_1}{\sup A_1-r_1}$, $\phi(t)=\left(1-\frac{\log(m)}{\log(m)+\log(1+\xi_0^2t)}\right)^2$,
$$ a(t)=e^{\alpha\phi(t)},\quad b(t)=e^{-\beta\phi(t)},\quad E(t,\varkappa)=1+(a(t)-1)\eta(\varkappa) ,$$
and
$$ A(t,\varkappa)=(1+\xi_0^2t)b(t)E(t,\varkappa),\quad B(t,\varkappa)=(1+\xi_0^2t)b(t)\cos(\varkappa) .$$

By construction, $\phi'(0)=0$, and hence $\frac{\partial g_{\mathscr{H};k,\xi_0,\epsilon}}{\partial t}\bigg|_{t=0}= 2\xi^2_0g_{\mathscr{H};k,\xi_0,\epsilon}$. By a straightforward calculation, we obtain
$$\phi'(t)=\frac{2\xi_0^2\log(m)\log(1+\xi_0^2 t)}{(1+\xi_0^2t)\left(\log(m)+\log(1+\xi_0^2t)\right)^3} ,$$
and 
$$\phi''(t)=\frac{2\xi_0^4\log(m)\big(\log(m)-2\log(1+\xi_0^2 t)-\log(m)\log(1+\xi_0^2 t)-(\log(1+\xi_0^2 t))^2\big)}{(1+\xi_0^2t)^2\left(\log(m)+\log(1+\xi_0^2t)\right)^4} .$$
Hence
\begin{eqnarray*}
\frac{1}{2}\frac{\partial g_{\mathscr{H};k,\xi_0,\epsilon}}{\partial t}\bigg|_{t=T}&\geq&\left(\frac{\xi_0^2}{1+\xi_0^2 T} +\frac{b'(T)}{b(T)}\right)g_{\mathscr{H};k,\xi_0,\epsilon}\\
&=&\frac{\xi_0^2}{1+\xi_0^2 T} \left(1 -\frac{2(1-\delta)\log(m)(\log(\frac{\xi_0^2}{r_1}))\log(1+\xi_0^2 T)}{\left(\log(m)+\log(1+\xi_0^2T)\right)^3}\right)g_{\mathscr{H};k,\xi_0,\epsilon}\\
&=&\frac{\xi_0^2}{1+\xi_0^2 T} \left(1 -\frac{100\delta^2\log(\frac{\xi_0^2}{r_1})}{\log(m)}\right)g_{\mathscr{H};k,\xi_0,\epsilon},
\end{eqnarray*}
where $T>0$ is the unique constant such that $\phi(T)=\frac{1}{1+\delta}$. Moreover, at $t=T$, we have
\begin{eqnarray*}
g_{\mathscr{H};k,\xi_0,\epsilon}\big|_{t=T}&=&(1+\xi_0^2 T)^2b^2(T)r_1^{-2}\cot^2(\xi_0)g_{\mathscr{A};k,\epsilon}|_{\partial\ball_{\xi_0}(x)} \\
& = & (1+\xi_0^2 T)^2\xi_0^{-4}\left(\xi_0^{\frac{8\delta}{1+\delta}} r_1^{\frac{-4\delta}{1+\delta}}\right)\cot^2(\xi_0)g_{\mathscr{A};k,\epsilon}|_{\partial\ball_{\xi_0}(x)}\\
&\geq& (1+\xi_0^2 T)^2\xi_0^{-4}\cot^2(\xi_0)g_{\mathscr{A};k,\epsilon}|_{\partial\ball_{\xi_0}(x)} ,
\end{eqnarray*}
and (d${}_2$) is proved.

It is therefore sufficient to prove (e${}_2$).

By Lemma \ref{lmmriccicurvature}, we can conclude that
\begin{eqnarray*}
\Ric_{g_{\mathscr{A};k,\epsilon}} &=& \left(1-\frac{f''(u)}{f(u)}\right)du^2+\left(\cos^2(u)+\frac{(k-4)f'(u)\sin(u)\cos(u)}{f(u)}\right) d\theta^2\\
&&+\left(k-5-(k-5)f'^2(u)-f(u)f''(u)+f(u)f'(u)\tan(u)\right)g_{\sphere^4} \\
&>&0.
\end{eqnarray*}

Now we estimate the Ricci curvature $\Ric_{g_{\mathscr{H};k,\xi_0,\epsilon}}$.

Let $v_0 =\frac{\partial }{\partial t} $, $v_1 =A^{-1} \frac{\partial }{\partial \varkappa} $, and $v_2 \in T_{\sphere^{k-4}} $ be unit vectors. By Lemma \ref{lmmriccicurvature3}, we can conclude that
\begin{eqnarray*}
\Ric_{g_{\mathscr{H};k,\xi_0,\epsilon}}(v_0,v_0)&=& -\frac{A_{tt}}{A}-(k-4)\frac{B_{tt}}{B}\\
&=&\left((k-3)\beta-\frac{\alpha\eta(\varkappa) a(t)}{E(t,\varkappa)}\right) \left(\phi''(t)+\frac{2\xi_0^2\phi'(t)}{1+\xi_0^2t}\right)\\
&&+\frac{\alpha\eta(\varkappa) a(t)\phi'(t)}{E(t,\varkappa)}\left(\frac{2\xi_0^2}{1+\xi_0^2t}-\alpha\phi'(t)\right)-(k-3)\beta^2\phi'^2(t),
\end{eqnarray*}
where the subscripts $t$ and $\varkappa$ denote partial differential with respect to $t$ and $\varkappa$, respectively. For example, $A_t = \frac{\partial A}{\partial t}$. By construction,
\begin{eqnarray*}
\phi''(t)+\frac{2\xi_0^2\phi'(t)}{1+\xi_0^2t}&=& \frac{2\xi_0^4\log(m)}{(1+\xi_0^2t)^2}\cdot \frac{\log(m)+(\log(m)-2)\log(1+\xi_0^2t)+(\log(1+\xi_0^2t))^2}{\left(\log(m)+\log(1+\xi_0^2t)\right)^4}\\
&>& \left(1-\frac{2}{\log(m)}\right)\cdot \frac{\xi_0^2\phi'(t)}{1+\xi_0^2t},
\end{eqnarray*}
\begin{eqnarray*}
\frac{2\xi_0^2}{1+\xi_0^2t}-\alpha\phi'(t)&=& \frac{2\xi_0^2}{1+\xi_0^2t}\cdot \frac{\left(\log(m)+\log(1+\xi_0^2t)\right)^3-\alpha\log(m)\log(1+\xi_0^2t)}{\left(\log(m)+\log(1+\xi_0^2t)\right)^3} \\
&\geq& (\log(m)-\alpha)\phi'(t)\geq 0,
\end{eqnarray*}
and
\begin{eqnarray*}
\beta^2\phi'(t)&=& \frac{2\beta^2\xi_0^2}{1+\xi_0^2t}\cdot \frac{\log(m)(\log(1+\xi_0^2t))}{\left(\log(m)+\log(1+\xi_0^2t)\right)^3} \leq\frac{2\xi_0^2}{1+\xi_0^2t}\cdot \frac{8\beta^2}{27\log(m)}.
\end{eqnarray*}
Since $\log(m)\geq \epsilon^{-10}\geq 100(k-3)\beta$ and $\beta\geq (1-\frac{10\log(\xi_0)}{\log(\epsilon)})\alpha$, we can conclude that
\begin{eqnarray*}
\Ric_{g_{\mathscr{H};k,\xi_0,\epsilon}}(v_0,v_0)&>&  \left((k-3)\beta-\alpha-\frac{(k-3)\beta-\alpha+1}{\log(m)}\right)\cdot \frac{\xi_0^2\phi'(t)}{1+\xi_0^2t}\\
&\geq& \frac{(k-3)\xi_0^2\phi'(t)}{2(1+\xi_0^2t)} \geq 0 .
\end{eqnarray*}

Similarly, by a straightforward calculation, we have $|\eta(\varkappa)\tan(\varkappa)|\leq \epsilon^{-4}$, and hence
\begin{eqnarray*}
|\Ric_{g_{\mathscr{H};k,\xi_0,\epsilon}}(v_0,v_1)|&=& (k-4)\frac{|A_tB_{\varkappa}-AB_{t{\varkappa}}|}{A^2B}\\
&=& \frac{(k-4)\alpha a(t)\phi'(t)|\eta(\varkappa)\tan(\varkappa)|}{(1+\xi_0^2t)b(t)E^2(t,\varkappa)}\\
&\leq& \frac{(k-4)\epsilon^{-4}\alpha a(t)\phi'(t)}{(1+\xi_0^2t)b(t)E^2(t,\varkappa)}.
\end{eqnarray*}

Now we consider $\Ric_{g_{\mathscr{H};k,\xi_0,\epsilon}}(v_1,v_1)$ and $\Ric_{g_{\mathscr{H};k,\xi_0,\epsilon}}(v_2,v_2)$.

Set $D(t,\varkappa)=E(t,\varkappa)-(a(t)-1)\eta'(\varkappa)\tan(\varkappa)$, 
$$U(t,\varkappa)=\frac{D(t,\varkappa)}{b^{2}(t)E^{3}(t,\varkappa)},\quad \textrm{and} \quad V(t,\varkappa)=\frac{E^2(t,\varkappa)-\sin^2(\varkappa)}{b^2(t)E^2(t,\varkappa)\cos^2(\varkappa)}.$$

Let $v_2'\in T_{\sphere^{k-4}}$ be another unit vector orthogonal to $v_2$. For each $t_0$, the sectional curvature on slides $\{t=t_0\}$ satisfies
$$ \sec_{t_0}(v_1,v_2)=\frac{U(t_0,\varkappa)}{(1+\xi_0^2t_0)^2},\quad \textrm{and} \quad \sec_{t_0}(v_2,v_2')=\frac{V(t_0,\varkappa)}{(1+\xi_0^2t_0)^2} .$$

For $V(t,\varkappa)$, by a straightforward calculation,
$$\frac{\partial\log V}{\partial t} = 2\beta\phi'(t)\left(1+\frac{\alpha a\eta\sin^2(\varkappa)}{\beta E(E^2-\sin^2(\varkappa))}\right)\geq 0,$$
and hence
$$V(t,\varkappa)\geq V(0,\varkappa)=1. $$

Now we consider $U(t,\varkappa)$. When $\eta'\tan(\varkappa)\leq 0$, we have $D(t,\varkappa)\geq E(t,\varkappa)$, and
\begin{eqnarray*}
U(t,\varkappa) &\geq & \frac{1}{b^2(t)E^2(t,\varkappa)} \geq \frac{1}{a^2(t)b^2(t)} = e^{2(\beta-\alpha)\phi(t)} \geq r_1^{20\delta} (\sup A_1)^{-2}\xi_0^4 \geq (1+\epsilon^{105})\xi_0^4 .
\end{eqnarray*}

Assume that $\eta'\tan(\varkappa)\geq 0$. By a direct computation,
\begin{eqnarray*}
\frac{\partial\log U}{\partial t}&=&\frac{\beta\phi'(t)}{E}\left(2(1-\eta)+2\eta a\left(1-\frac{\alpha}{\beta}\right)-\frac{\alpha a\eta'\tan(\varkappa)}{\beta D}\right)\\
&=&\beta\phi'(t)\left(2-\frac{\alpha a}{\beta E}\left(2\eta+\frac{\eta'\tan(\varkappa)}{D}\right)\right),
\end{eqnarray*}
and
\begin{eqnarray*}
&&\frac{\partial}{\partial t}\left(2(1-\eta)+2\left(1-\frac{\alpha}{\beta}\right)\eta a-\frac{\alpha a\eta'\tan(\varkappa)}{\beta D}\right)\\
&=&\alpha \phi'(t)a\left(2\left(1-\frac{\alpha}{\beta}\right)\eta-\frac{\alpha\eta'\tan(\varkappa)(1-\eta+\eta'\tan(\varkappa))}{\beta D^2}\right)\leq 0 .
\end{eqnarray*}
Then $\frac{\partial}{\partial t}\left(\frac{E}{\beta\phi'(t)}\cdot\frac{\partial\log U}{\partial t}\right)\leq 0 $. Hence 
\begin{eqnarray*}
\inf_{t\in[0,T]} U(t,\varkappa) &=& \min\{U(0,\varkappa),U(T,\varkappa)\}\geq \min\left\{1,(1+\epsilon^{100})r_1^2b^{-2}(T)\right\}\\
&=& \min\left\{1,(1+\epsilon^{100})(r_1\xi_0^{-2})^{\frac{4\delta}{1+\delta}}\xi_0^4\right\} \geq (1+\epsilon^{100}) \xi_0^4 .
\end{eqnarray*}
Here we used $\sec\left( g_{\mathscr{A};k,\epsilon}|_{\partial\ball_{\xi_0}(x)}\right)\geq (1+\epsilon^{100})\cot^2 (\xi_0)$. 

Since $|\phi'|\leq \frac{16\xi_0^2}{27(1+\xi_0^2 t)\log(m)}$ and $|\phi''|\leq \frac{10\xi_0^4}{27(\log(m))2}$, we can conclude that
\begin{eqnarray*}
\Ric_{g_{\mathscr{H};k,\xi_0,\epsilon}}(v_1,v_1) &=&\frac{2\beta\xi_0^2\phi'}{1+\xi_0^2 t}-\beta^2\phi'^2(t)+\beta\phi''(t)-\frac{\alpha a\eta}{E}
\left(\phi''(t)+\alpha\phi'^2(t)\right)\\
&&+\frac{(k-4)(U-\xi_0^4)}{(1+\xi_0^2 t)^2}+\frac{(k-4)\xi_0^2\phi'(t)}{1+\xi_0^2 t}\left(2\beta-\frac{\alpha a\phi'(t)\eta}{E}\right)\\
&\geq & \frac{\epsilon^{200}\xi_0^4}{(1+\xi_0^2 t)^2} .
\end{eqnarray*}
It follows that
\begin{eqnarray*}
&&\Ric_{g_{\mathscr{H};k,\xi_0,\epsilon}}(v_0,v_0)\Ric_{g_{\mathscr{H};k,\xi_0,\epsilon}}(v_1,v_1)-\Ric_{g_{\mathscr{H};k,\xi_0,\epsilon}}(v_0,v_1)^2\\
&\geq & \frac{(k-3)\xi_0^2\phi'(t)}{2(1+\xi_0^2t)} \cdot \frac{\epsilon^{200}\xi_0^4}{(1+\xi_0^2 t)^2} - \frac{(k-4)^2\epsilon^{-8}\alpha^2 a^2(t)\phi'^2(t)}{(1+\xi_0^2t)^2 b^2(t)E^4(t,\varkappa)} \\
&\geq & \frac{(k-4)\phi'(t)}{(1+\xi_0^2t)^2} \cdot \left( \frac{\epsilon^{200}\xi_0^6}{2(1+\xi_0^2 t)} - \frac{(k-4)\epsilon^{-8}\alpha^2 a^2(t)\phi'(t)}{b^2(t)E^4(t,\varkappa)} \right) \\
&\geq & \frac{(k-4)\phi'(t)}{(1+\xi_0^2t)^2} \cdot \left( \frac{\epsilon^{200}\xi_0^6}{2(1+\xi_0^2 t)} - (k-4)\epsilon^{-8}\alpha^2 e^{2(\alpha+\beta)}\phi'(t) \right) >0.
\end{eqnarray*}

Similarly, we have
\begin{eqnarray*}
\Ric_{g_{\mathscr{H};k,\xi_0,\epsilon}}(v_2,v_2)&=&\beta\left(\phi''(t)+\frac{2\xi_0^2\phi'(t)}{1+\xi_0^2 t}\right) -\beta^2\phi'^2(t) +\frac{U-\xi_0^4 +(k-5)(V-\xi_0^4)}{(1+\xi_0^2 t)^2} \\
& & + \left(\frac{2(k-4)\xi_0^2\beta-\xi_0^2\alpha a\eta}{1+\xi_0^2 t}+\frac{\alpha\beta a\eta \phi'(t)}{E}-(k-4)\beta^2\phi'(t)\right)\phi'(t) \\
&\geq & \frac{\epsilon^{200}}{(1+\xi_0^2 t)^2} >0.
\end{eqnarray*}

By the warped product structure, we have
$$\Ric_{g_{\mathscr{H};k,\xi_0,\epsilon}}(v_0,v_2)=\Ric_{g_{\mathscr{H};k,\xi_0,\epsilon}}(v_1,v_2) =0,$$
and hence $\Ric_{g_{\mathscr{H};k,\xi_0,\epsilon}}>0$. This completes the proof.
\end{proof}

The following proposition is the key step in the proof of Proposition \ref{prop_dockingstation}. Here we use the method of constructing hypersurfaces via involutions from \cite{mrs1} to avoid the smoothing of corners, thereby simplifying the construction in \cite{burdick1}.

\begin{prop}
\label{Prop_PerelmanConstruction_threeholes}
For any $k\geq 6$ and $\rho\in(0,1)$, there exist a $(k-2)$-dimensional compact connected smooth Riemannian manifold $(W_{k,\rho},g_{W_{k,\rho}})$ with boundary $\partial{W}_{k,\rho}$ and an involution $\sigma_{k,\rho}:(W_{k,\rho},g_{W_{k,\rho}})\to (W_{k,\rho},g_{W_{k,\rho}})$, such that:
\begin{enumerate}[\rm (a${}_3$).]
\item The Ricci curvature $\Ric_{g_{W_{k,\rho}}}>0$,
\item The manifold $W_{k,\rho}$ is diffeomorphic to $\sphere^{k-3}$ with three pairwise disjoint open balls removed, and the boundary $\partial W_{k,\rho}$ has three connected components $V_{k,\rho;1}\sqcup V_{k,\rho;2}\sqcup V_{k,\rho;3}$,
\item For $j=1,2,3 $, $(V_{k,\rho;j},g_{g_{W_{k,\rho}}}|_{V_{k,\rho;j}})$ is isometric to $\sphere_{r_{k,\rho}}^{k-3}$ for some $r_{k,\rho}>0$, and its outer second fundamental form $\mathbf{II}_{k,\rho;j}\geq -\rho r_{k,\rho}^{-1} g_{W_{k,\rho}}|_{V_{k,\rho;j}}$,
\item The restriction $\sigma_{k,\rho}|_{V_{k,\rho;1}}$ is the involution that fixes the equator,
\item The restriction $\sigma_{k,\rho}|_{V_{k,\rho;2}}$ gives an isomorphism $V_{k,\rho;2}\to V_{k,\rho;3}$,
\item The set of fixed points $H_{k,\rho}=\{x\in W_{k,\rho}:\sigma_{k,\rho}(x)=x\}$ is a hypersurface in $W_{k,\rho}$,
\item The metric $g_{W_{k,\rho}}|_{H_{k,\rho}}$ can be represented by
$$g_{W_{k,\rho}}|_{H_{k,\rho}} =dr_*^2+\varphi_{H_{k,\rho}}^2(r_*),\quad r_*\in [0,L_{k,\rho}] ,$$
where $\varphi_{H_{k,\rho}}\geq 0$, $\varphi_{H_{k,\rho}}|_{[0,L_{k,\rho})}>0$, $\varphi_{H_{k,\rho}}''|_{[0,L_{k,\rho})}<0$, $\varphi_{H_{k,\rho}} (0)=r_{k,\rho}$, $\varphi_{H_{k,\rho}}'(0)\leq \rho$, and $r_*(x)=\dist (V_{k,\rho;1},x)$ on a small neighborhood of $V_{k,\rho;1}$.
\end{enumerate}
\end{prop}

\begin{proof}
At first, we describe the construction of $(W_{k,\rho},g_{W_{k,\rho}})$ and $\sigma_{k,\rho}$.

Fix $\xi_0=e^{-10k}\rho^{10}$ and $\epsilon=e^{-200k}\xi_0^{200}$. Let $g_{\mathscr{A};k,\epsilon}$ and $g_{\mathscr{H};k,\xi_0,\epsilon}$ denote the metrics on $\sphere^{k-2}$ and $[0,T_{\mathscr{H};k,\xi_0,\epsilon}]\times \sphere^{k-3}$, respectively, constructed in Lemma \ref{lmm_Perelman_ambientandcollar}. Recall that $g_{\mathscr{A};k,\epsilon}$ can be represented by
$$g_{\mathscr{A};k,\epsilon}=du^2+\cos^2 (u)d\theta^2+f^2(u) g_{\sphere^{k-4}},$$
where $u\in [0,\frac{\pi}{2}]$, $\theta\in\mathbb{R}/2\pi\mathbb{Z}$, $f|_{(0,\frac{\pi}{2}]}>0$, $f'|_{[0,\frac{\pi}{2})}>0$, $f''|_{(0,\frac{\pi}{2}]}<0$, $f'(0)=1$, $f(\frac{\pi}{2})=\epsilon$, and $f(0)=f^{(\mathrm{even})}(0)=f^{(\mathrm{odd})}(\frac{\pi}{2})=0$. By the $\sphere^1$-action given by $\theta\in\mathbb{R}/2\pi\mathbb{Z}$, the involution $\theta\mapsto-\theta$ gives an involution $\sigma_{k,\epsilon;1}:(\sphere^{k-2},g_{\mathscr{A};k,\epsilon})\to (\sphere^{k-2},g_{\mathscr{A};k,\epsilon})$, and the set of fixed points is $\{\theta=0\}\cup\{\theta=\pi\}\cup\{u=\frac{\pi}{2}\}$. Let $x_1\in\sphere^{k-2}$ such that $u(x_1)=\theta(x_1)=0$, and let $x_2,x_3\in\sphere^{k-2}$ such that $u(x_2)=u(x_3)=0$ and $\theta(x_2)=-\theta(x_3)=\frac{1}{10}$.

By construction, the Ricci curvature of $g_{\mathscr{A};k,\epsilon}$ and $g_{\mathscr{H};k,\xi_0,\epsilon}$ are strictly positive, the metric $g_{\mathscr{H};k,\xi_0,\epsilon}\big|_{t=T_{\mathscr{H};k,\xi_0,\epsilon}}$  is isometric to $R_{\mathscr{H};k,\xi_0,\epsilon}^2g_{\mathscr{A};k,\epsilon}\big|_{\partial \ball_{\xi_0}(x_j)}$, $j=1,2,3$, and the outer second forms satisfies $\mathbf{II}_{\{t=T_{\mathscr{H};k,\xi_0,\epsilon}\}}>R_{\mathscr{H};k,\xi_0,\epsilon}^{-1}\mathbf{II}_{\partial \ball_{\xi_0}(x_j)}$. By gluing $[0,T_{\mathscr{H};k,\xi_0,\epsilon}]\times \sphere^{k-3}$ to $\sphere^{k-2}\setminus \cup_{j=1}^3 \ball_{\xi_0}(x_j)$ along the boundary, one can get a smooth manifold $W_{k,\rho}$ with boundary $\partial W_{k,\rho}$. Moreover, the glued metric $\hat{g}_{W_{k,\rho}}$ is a $C^{0,1}$ metric. By Lemma \ref{lmmgluing}, we can modify $\hat{g}_{W_{k,\rho}}$ near the glued hypersurface and obtain a smooth Riemannian metric $g_{W_{k,\rho}}$ with positive Ricci curvature. This is the construction of $W_{k,\rho}$. By the construction, we can extend $\sigma_{k,\epsilon;1}|_{sphere^{k-2}\setminus \cup_{j=1}^3 \ball_{\xi_0}(x_j)}$ to an involution $\sigma_{k,\rho}:W_{k,\rho}\to W_{k,\rho}$.

Properties (a${}_3$)-(f${}_3$) follow from the construction. It suffices to prove (g${}_3$) now. By construction, the metric $\hat{g}_{W_{k,\rho}}|_{H_{k,\rho}}$ can be represented by
$$\hat{g}_{W_{k,\rho}}|_{H_{k,\rho}} =dr_*^2+\hat{\varphi}_{H_{k,\rho}}^2(r_*),\quad r_*\in [0,L_{k,\rho}] ,$$
where $\hat{\varphi}$ is smooth on $[0,L_{k,\rho}]\setminus\{\hat{t}_{k,\rho}\}$ for some $\hat{t}_{k,\rho}\in (0,L_{k,\rho})$, $\hat{\varphi}_{H_{k,\rho}}|_{[0,L_{k,\rho})}>0$, $\hat{\varphi}_{H_{k,\rho}} (0)=r_{k,\rho}$, $\hat{\varphi}_{H_{k,\rho}}'(0)\leq \rho$, $\hat{\varphi}_{H_{k,\rho}}'(\hat{t}_{k,\rho}^{+})<\hat{\varphi}_{H_{k,\rho}}'(\hat{t}_{k,\rho}^{-})$, $\hat{\varphi}_{H_{k,\rho}}''|_{(\hat{t}_{k,\rho},L_{k,\rho})}<0$, and $r_*(x)=\dist (V_{k,\rho;1},x)$ on a small neighborhood of $V_{k,\rho;1}$. Moreover, we have $\hat{\varphi}_{H_{k,\rho}}|_{(0,\hat{t}_{k,\rho})}=\kappa'(1+\xi_0^2 t)b(t)$ for some $\kappa'>0$, where
$$b(t)=e^{-\beta\phi(t)},\quad \beta=(1-\delta)(\log(\xi_0^2)-\log (r_1)) ,\quad \phi(t)=\left(1-\frac{\log(m)}{\log(m)+\log(1+\xi_0^2t)}\right)^2,$$
$t=\kappa'' r_*$ for some $\kappa''>0$, and $r_1,m>0$ are constants defined in Lemma \ref{lmm_Perelman_ambientandcollar}. 

By a straightforward calculation,
\begin{eqnarray*}
\big((1+\xi_0^2 t)b\big)''(t) &=& 2\xi_0^2 b'(t)+(1+\xi_0^2 t)b''(t) \\
&=& -\beta (1+\xi_0^2 t)b(t) \left(\phi''(t)+\frac{2\beta\xi_0^2\phi'(t)}{1+\xi_0^2 t}-\beta\phi'^2(t) \right) <0,
\end{eqnarray*}
and hence $\hat{\varphi}''_{H_{k,\rho}}|_{(0,\hat{t}_{k,\rho})}<0$. Now we can apply the construction in Lemma \ref{lmm_Perelman_ambientandcollar} to show that $\varphi''<0$. This proves (g${}_3$), and the proposition follows.
\end{proof}

Now we are ready to prove Proposition \ref{prop_dockingstation}.

\vspace{0.2cm}

\noindent \textbf{Proof of Proposition \ref{prop_dockingstation}: }
Fix $\rho=e^{-100k^{100}}\epsilon_0^{100}$. Let $W_{k,\rho}$, $H_{k,\rho}$, $\sigma_{k,\rho}$ and $\{(V_{k,\rho;j},g_{g_{W_{k,\rho}}}|_{V_{k,\rho;j}})\}_{j=1}^3$ be as in Proposition \ref{Prop_PerelmanConstruction_threeholes}. Let 
$$\{g_{2,k-3;i}\}_{i=1}^\infty=\left\{dr^2+\varphi_{2,k-3;i}^2(r)g_{\sphere^2}+\psi_{2,k-3;i}^2(r)g_{\sphere^{k-3}}\right\}_{i=1}^\infty$$ 
be the sequence of Riemannian metrics on $\mathbb{R}^{k-2}$ defined in Lemma \ref{lmm_firstdoubleproduct}.

After scaling, one can glue $\mathfrak{B}_{\frac{\pi}{4};i}=\{r\leq \frac{\pi}{4}\}\subset (\mathbb{R}^3\times\sphere^{k-3},g_{2,k-3;i})$ to $W_{k,\rho}\times \sphere_{r_{k,\rho}}^2$ along $V_{k,\rho;1}\times \sphere_{r_{k,\rho}}^2$, and obtain a smooth manifold $\mathfrak{W}_{k,\rho}$ equipped with a $C^{0,1}$ metric. By applying Lemma \ref{lmmgluing} to the glued metric, one can obtain a smooth metric $\tilde{g}_{k,\rho}$ on $\tilde{W}_{k,\rho}$. Moreover, the involution $\sigma_{k,\rho}$ can be extended to $\tilde{\sigma}_{k,\rho}:\mathfrak{W}_{k,\rho}\to\mathfrak{W}_{k,\rho}$, and the set of fixed points of $\tilde{\sigma}_{k,\rho}$ is a hypersurface $\mathfrak{H}_{k,\rho}$ in $\mathfrak{W}_{k,\rho}$ diffeomorphic to $\sphere^{k-1}$. Clearly, $\mathfrak{H}_{k,\rho}$ devides $\mathfrak{W}_{k,\rho}$ into two components. Let $\mathfrak{Y}$ be one of the components.

By construction, $(\mathfrak{Y},\tilde{g}_{k,\rho})$ is a smooth Riemannian manifold with boundary $\mathfrak{Z}_1\cup\mathfrak{Z}_2$. By scaling, we can assume that $\mathfrak{Z}_1$ is isometric to $\sphere_1^{k-3}\times\sphere_1^2$, and its outer second fundamental form $\mathbf{II}_{\mathfrak{Z}_1}\geq -\epsilon_0 \tilde{g}_{k,\rho}$. Moreover, $\mathfrak{Z}_2=\mathfrak{H}_{k,\rho}$ is diffeomorphic to $\sphere_1^{k-1}$, the Ricci curvature $\Ric_{\tilde{g}_{k,\rho}|_{\mathfrak{Z}_2}}>0 $, and the outer second fundamental form $\mathbf{II}_{\mathfrak{Z}_2}=0$. By Lemma \ref{lmm_smallpertubation_positivesecond}, we can modify the metric near $\mathfrak{Z}_2$ such that the outer second fundamental form $\mathbf{II}_{\mathfrak{Z}_2}$ becomes strictly positive.

By Lemma \ref{lmmcoldingnaberexamplelemma}, we only need to find a continuous family of smooth metrics $\tilde{g}(t)$ on $\sphere_1^{k-1}$ with positive Ricci curvature, $t\in[0,2]$, such that $\tilde{g}(0)\cong\tilde{g}_{k,\rho}|_{\mathfrak{Z}_2}$, and $\tilde{g}(2)$ is the round metric.

By construction, after rescaling, the metric $\tilde{g}_{k,\rho}|_{\mathfrak{Z}_2}$ can be represented by
$$\tilde{g}_{k,\rho}|_{\mathfrak{Z}_2} =dr^2 + \varphi (r)^2 g_{\sphere^2} +  \psi (r)^2 g_{\sphere^{k-4}},\quad r\in \left[0,\frac{\pi}{2}\right] ,$$
where $\varphi(0)=\varphi^{(\mathrm{even})}(0)=\psi(\frac{\pi}{2})=\psi^{(\mathrm{even})}(\frac{\pi}{2})=0$, $\varphi^{(\mathrm{odd})}(\frac{\pi}{2})=\psi^{(\mathrm{odd})}(0)=0$, $\varphi'(0)=-\psi'(\frac{\pi}{2})=1$, $\psi'\in [-1,0]$, $\lim_{r\to 0^+}\frac{\varphi''}{\varphi}(r)<0$, $\varphi''\leq 0$, $\varphi(r)\geq r|\varphi'(r)|$, $\psi(r)\geq (\frac{\pi}{2}-r)|\psi'(r)|$, $|\varphi|+|\psi|<\rho$ and $\Ric_{\tilde{g}_{k,\rho}|_{\mathfrak{Z}_2}}>0$. Moreover, there exists a constant $L\in(0,\frac{\pi}{2})$ such that the restriction $\psi''|_{(L,\frac{\pi}{2}]}< 0$, $\psi'|_{(L,\frac{\pi}{2}]}\leq 0$, and for any $r\in[0,L]$, $|\psi'(r)|+\psi(r)|\psi''(r)|\leq\rho$. Without loss of generality, we can assume that $\varphi_t (r)\leq \rho^{10}\cdot\inf_{[0,L]}\psi $, $\forall r\in[0,L]$. 

Let $\tilde{g}(t)=dr^2 + \varphi_t^2 (r) g_{\sphere^2} + \psi_t^2 (r) g_{\sphere^{k-4}}$, $t\in[0,1]$, where
$$\varphi_t (r)=\varphi (r),\quad \psi_t (r)=t\varphi \left(\frac{\pi}{2}-r\right)+(1-t)\psi(r).$$
Set $X_0=\frac{\partial}{\partial r}$, $X_1\in T_{\sphere^2}$ and $X_2\in T_{\sphere^{k-4}}$. By Lemma \ref{lmmriccicurvature}, the Ricci curvature tensor can be determined by
\begin{eqnarray*}
\Ric_{\tilde{g}(t)} \left(X_0 \right) & = & - \left( \frac{2\varphi_t''(r)}{\varphi_t(r)} + (k-4)\frac{\psi_t'' (r)}{\psi_t (r)} \right) X_0 ,\\
\Ric_{\tilde{g}(t)} \left(X_1 \right) & = & \left(-  \frac{\varphi_t''(r)}{\varphi_t(r)} + \frac{1-\varphi_t'^2(r) }{\varphi_t^2(r)} - (k-4) \frac{\varphi_t'(r)\psi_t' (r))}{\varphi_t(r)\psi_t (r) } \right) X_1 ,\\
\Ric_{\tilde{g}(t)} \left(X_2 \right) & = & \left(- \frac{\psi_t'' (r)}{\psi_t (r)} + (k-5)\frac{1-\psi_t'^2 (r) }{\psi_t^2(r)} - \frac{2\varphi_t'(r)\psi_t' (r))}{\varphi_t(r)\psi_t (r) } \right) X_2 .
\end{eqnarray*}
By estimating in the cases of $r\in[0,L]$ and $r\in[L,\frac{\pi}{2}]$ respectively, one can obtain $\Ric_{\tilde{g}(t)}>0$, $\forall t\in[0,1]$. For any $t\in[1,2]$, let $\tilde{g}(t)=dr^2 + \varphi_t^2 (r) g_{\sphere^2} + \psi_t^2 (r) g_{\sphere^{k-4}}$, $t\in[1,2]$, where
$$\varphi_t (r)=(t-1)\sin (r)+(2-t)\varphi (r),\quad \psi_t (r)=(t-1)\cos (r)+(2-t)\varphi \left(\frac{\pi}{2}-r\right) .$$
Clearly, $\Ric_{\tilde{g}(t)}>0$, $\forall t\in[1,2]$. This completes the proof.
\qed

\subsection{Warped product structure}
In the proof of Proposition \ref{prop_eightdim_localmodel}, the key step is to construct suitable Riemannian metrics with positive Ricci curvature on the first component $\mathbb{R}^k$ as follows.

\begin{prop}
\label{prop_sixdim_localmodel}
Let $ k\geq 6 $. Then there exist:
\begin{enumerate}[\rm (i${}_{2}$).]
\item Three precompact open subset $U\subset\bar{U}\subset V\subset\bar{V}\subset W\subset\bar{W}\subset \mathbb{R}^k$ with smooth boundaries,
\item A sequence of complete Riemannian manifolds $\{h_j\}_{j=1}^\infty$ on $\mathbb{R}^k$,
\item A sequence of pairs of functions $\{\varphi_j,\psi_j\}_{j=1}^\infty$ satisfying (\ref{conditionvarphi}) and (\ref{conditionphirho}),
\end{enumerate}
such that
\begin{enumerate}[\rm (i${}_{2}$).]
\addtocounter{enumi}{4}
\item For any $\kappa\in(0,\frac{1}{10})$, the functions $\varphi_j(r)$ and $\psi_j(r)$ converge to $\rho\sin (r)$ smoothly on $[\kappa,\frac{\pi}{5}]$ as $j\to\infty$, where $\rho=\rho(k,\epsilon,\alpha,\delta)>0$ is a constant,
\item The restriction of $h_j$ on $V$ can be represented by
$$\big(V,h_j|_{V}\big)\cong\left([0,\frac{\pi}{6})\times \sphere^{k-3}\times\sphere^2,dr^2+\varphi_i^2(r)g_{\sphere^{k-3}}+\psi_i^2(r)g_{\sphere^2}\right) ,$$
\item The metrics $h_j$ converge to a Riemannian metric $h_\infty$ on $\mathbb{R}^k\setminus U$ smoothly,
\item The restriction of $h_j$ on $\mathbb{R}^k\setminus W$ is independent of $j$, and it can be represented by
$$(\mathbb{R}^k\setminus W,h_j|_{\mathbb{R}^k\setminus W})\cong\left([R,\infty)\times \sphere^{k-1},dr^2+\epsilon^2(r+\arctan(r))^2g_{\sphere^{k-1}}\right) ,$$
where $\epsilon\in(0,1)$ and $R>1$ are positive constants,
\item The Ricci curvature $\Ric_{h_{j}}\geq \frac{\tau}{1+r_j^4}$ for some constant $\tau>0$ independent of $j$, where $r_j(x)=\dist_{h_j}(0^k,x)$.
\end{enumerate}
\end{prop}

\begin{proof}
For any small constant $\rho\in(0,e^{-100e^{100k}})$, let $(\hat{Y}_{k,\rho},g_{\hat{Y}_{k,\rho}})$ be the Riemannian manifold with boundary $\partial{\hat{Y}}_{k,\rho}$ constructed in Proposition \ref{prop_dockingstation}, such that:
\begin{itemize}
\item The Ricci curvature $\Ric_{g_{\hat{Y}_{k,\rho}}}>0$,
\item The manifold $\hat{Y}_{k,\rho}$ is diffeomorphic to $\ball_2^k(0)\setminus \ball_{e^{-10k}} (\sphere_1^{2}) $, where $\ball_2^k(0)$ is the ball centered at $0\in\mathbb{R}^k$ with radius $2$, and $\sphere_1^{2}\subset\mathbb{R}^3\times\{0\}\subset \mathbb{R}^k$ is the unit sphere.
\item The boundary $\partial\hat{Y}_{k,\rho}$ has two connected components $\hat{Z}_{k,\rho;1}\sqcup \hat{Z}_{k,\rho;2}$,
\item $(\hat{Z}_{k,\rho;1},g_{\hat{Y}_{k,\rho}}|_{\hat{Z}_{k,\rho;1}})$ is isometric to $\sphere_1^{k-3}\times\sphere_1^2$ and its outer second fundamental form $\mathbf{II}_{k,\rho;1}\geq -\rho g_{\hat{Y}_{k,\rho}}|_{\hat{Z}_{k,\rho;1}}$,
\item $(\hat{Z}_{k,\rho;2},g_{\hat{Y}_{k,\rho}}|_{\hat{Z}_{k,\rho;2}})$ is isometric to $\sphere^{k-1}_{R_{k,\rho}}$ and its outer second fundamental form $\mathbf{II}_{k,\rho;2}\geq \delta g_{\hat{Y}_{k,\rho}}|_{\hat{Z}_{k,\rho;2}}$, where $R_{k,\rho}>1$, $\delta\in(0,\rho)$ are constants.
\end{itemize}

Let $\{g_{k-3,2;i}\}_{i=1}^\infty$ be the sequence of complete Riemannian metrics constructed in Lemma \ref{lmm_firstdoubleproduct}. Then $g_{k-3,2;i}$ can be represented by
$$g_{k-3,2;i}=dr^2+\varphi_{k-3,2;i}^2(r)g_{\mathbb{S}^{k-3}}+\psi_{k-3,2;i}^2(r) g_{\mathbb{S}^2}$$
on $\mathbb{R}^{k-2}\times\mathbb{S}^2$, where $\varphi_{k-3,2;i},\psi_{k-3,2;i}$ are nonnegative increasing functions satisfying (\ref{conditionvarphi}) and (\ref{conditionphirho}), such that
$$ \varphi_{k-3,2;i}(r)=\psi_{k-3,2;i}(r)=e^{-10k}\sin(r+\epsilon_i) \quad\textrm{ for } r\in \left[e^{-10(k+i)},\frac{\pi}{4}\right],$$
and the Ricci curvature $\Ric_{g_{k-3,2;i}} \geq\frac{\tau_0}{1+r^2} $, where $\tau_0>0$ is a positive constant, and $\{\epsilon_i\}_{i=1}^\infty$ is a sequence of positive constants converging to $0$.

Now we glue $\left(\{r\leq \frac{\pi}{4}-\epsilon_i\},\sqrt{2}e^{10k}g_{k-3,2;i}\right)$ to $(\hat{Y}_{k,\rho},g_{\hat{Y}_{k,\rho}})$ along $\hat{Z}_{k,\rho;1}$, and glue
$$\left([\delta^{-10}R_{k,\rho},\infty)\times \sphere^{k-1},dr^2+\epsilon(r+\arctan(r))^2g_{\sphere^{k-1}}\right) $$
to $(\hat{Y}_{k,\rho},g_{\hat{Y}_{k,\rho}})$ along $\hat{Z}_{k,\rho;2}$, where $\epsilon=R_{k,\rho}(\delta^{-10}R_{k,\rho}+\arctan(\delta^{-10}R_{k,\rho}))^{-1}$. By Lemma \ref{lmmgluing}, after modifying the metric near $\partial\hat{Y}_{k,\rho}$, one can obtain a smooth Riemannian metric with positive Ricci curvature.

Let $U=\{r\leq \frac{\pi}{8}\}\subset\left(\mathbb{R}^{k-2}\times\sphere^2,\sqrt{2}e^{10k}g_{k-3,2;i}\right)$, and let $W\subset \mathbb{R}^k$ containing $\hat{Y}_{k,\rho}$ such that 
$$\mathbb{R}^k\setminus W\cong \left([2\delta^{-10}R_{k,\rho},\infty)\times \sphere^{k-1},dr^2+\epsilon(r+\arctan(r))^2g_{\sphere^{k-1}}\right).$$
All the properties now follow from Lemma \ref{lmm_firstdoubleproduct} and Proposition \ref{prop_dockingstation}.
\end{proof}

By combining Lemma \ref{lmm3dimlocalmodel} and Proposition \ref{prop_sixdim_localmodel}, one can prove Proposition \ref{prop_eightdim_localmodel}. For convenience, we recall the statement of Proposition \ref{prop_eightdim_localmodel} here.

\begin{prop}[Imply Proposition \ref{prop_eightdim_localmodel}]
\label{prop_eightdim_localmodel2}
Let $ k\geq 6 $. Then for any $0< \epsilon \leq \epsilon_0(k) $, $0<\alpha \leq \alpha_0 (k,\epsilon )<1 $ and $0<\delta< \delta_0 (k,\epsilon,\alpha )$, there exist:
\begin{enumerate}[\rm (i${}_{3}$).]
\item Three precompact open subset $U\subset\bar{U}\subset V\subset\bar{V}\subset W\subset\bar{W}\subset \mathbb{R}^k$ with smooth boundaries,
\item A sequence of warped product metrics $\{g_i=h_i+f^2g_{\sphere^{2}}\}_{j=1}^\infty$, where $\{h_i\}_{i=1}^\infty$ is a sequence of complete Riemannian metrics on $\mathbb{R}^k$, $f$ is a smooth positive function on $\mathbb{R}^k$, and $g_{\sphere^2}$ denotes the round metric on $\sphere^2$,
\item A sequence of pairs of functions $\{(\varphi_i,\psi_i)\}_{i=1}^\infty$ satisfying (\ref{conditionvarphi}) and (\ref{conditionphirho}),
\end{enumerate}
such that
\begin{enumerate}[\rm (i${}_{3}$).]
\addtocounter{enumi}{3}
\item For any $\kappa\in(0,\frac{1}{10})$, the functions $\varphi_i(r)$ and $\psi_i(r)$ converge to $\rho\sin (r)$ smoothly on $[\kappa,\frac{\pi}{8}]$ as $i\to\infty$, where $\rho=\rho(k,\epsilon,\alpha,\delta)>0$ is a constant,
\item The metrics $h_i$ converge to a Riemannian metric $h_\infty$ on $\mathbb{R}^k\setminus U$ smoothly,
\item The restriction of $h_i$ on $V$ can be represented by
$$\big(V,h_i|_{V}\big)\cong\left([0,\frac{\pi}{10})\times \sphere^{k-3}\times\sphere^2,dr^2+\varphi_i^2(r)g_{\sphere^{k-3}}+\psi_i^2(r)g_{\sphere^2}\right) ,$$
\item The restriction of $g_i$ on $(\mathbb{R}^k\setminus W)\times \sphere^2$ is independent of $i$, and it can be represented by
$$\big((\mathbb{R}^k\setminus W)\times \sphere^2,g_i|_{(\mathbb{R}^k\setminus W)\times \sphere^2}\big)\cong\left([R,\infty)\times \sphere^{k-1}\times\sphere^2,dr^2+(1-\epsilon)^2r^2g_{\sphere^{k-1}}+\delta^2r^{2\alpha}g_{\sphere^2}\right) ,$$
where $R=R(k,\epsilon,\alpha,\delta)>1$ is a positive constant,
\item The Ricci curvature $\Ric_{g_{i}}\geq \frac{\tau}{1+r_i^2}$ for some constant $\tau=\tau(k,\epsilon,\alpha,\delta)>0$, where $r_j$ is the distance function defined by $r_i(x)=\dist_{g_i}(\{0^k\}\times\sphere^2,x)$.
\end{enumerate}
\end{prop}

\begin{proof}
Let $U_0,V_0,W_0$ and $h_{i,0}$ be the data defined in Proposition \ref{prop_sixdim_localmodel}. Set
$$g_{i,0}=h_{i,0}+g_{\sphere^2},\qquad a(r)=\epsilon_0(r+\arctan r),$$
where $\epsilon_0$ is the positive constant in the infinite annulus of Proposition \ref{prop_sixdim_localmodel}. Then, on $\mathbb{R}^k\setminus W_0$, the metric $g_{i,0}$ is independent of $i$ and has the form
\begin{equation*}
g_{i,0}=dr^2+a(r)^2g_{\sphere^{k-1}}+g_{\sphere^2},\quad r\in[R_0,\infty),
\end{equation*}
where $R_0>1$ is a constant.

We now apply Lemma \ref{lmm3dimlocalmodel} with $p=k-1$ and $q=2$. For the given $\epsilon,\alpha,\delta$, this gives functions $\varphi,\psi$ and a metric
$$\hat{g}=dt^2+\varphi^2(t)g_{\sphere^{k-1}}+\psi(s)^2g_{\sphere^2}$$
such that
$$\Ric_{\hat{g}}\geq \frac{\tau_1}{1+t^2},\quad \varphi(0)=a(2R_0),\quad \psi(0)=1,\quad \varphi'(0)=a'(2R_0),\quad \psi'(0)=0,$$
and for $t\geq R_1$,
$$\varphi(t)=(1-\epsilon)(r+C),\quad \psi(t)=\delta (t+C)^\alpha,$$
where $\tau_1\in(0,1)$, $R_1>1$, and $C>0$ are constants.

Along the boundary $\{r=2R_0\}\subset\mathbb{R}^k\setminus W_0$, we glue $g_{i,0}$ to $\hat{g}$ by Lemma \ref{lmm_smallpertubation_positivesecond} and Lemma \ref{lmmgluing}. The The modification is supported in a small neighborhood of $\{r=2R_0\}$, and hence is independent of $i$.

Let $g_i$ denote the glued metric. Outside the gluing domain it coincides either with $g_{i,0}$ or with $\hat{g}$. Moreover, it can be represented by
$$g_i=h_i+f^2g_{\sphere^2},$$
where $f=1$ on the part unchanged from $g_{i,0}$, and $f=\psi$ on the doubly warped product end. Similarly, on the doubly warped product end, we have
$$h_i=dt^2+\varphi^2(t)g_{\sphere^{k-1}}.$$

Since all modifications are on $\mathbb{R}^k\setminus W_0$, the metrics $h_i=h_{i,0}$ on $V_0$. The convergence property therefore follows from Proposition \ref{prop_sixdim_localmodel}. Similarly, $h_i$ converges smoothly on $\mathbb{R}^k\setminus U_0$.

It remains to prove the Ricci lower bound. Since the Ricci curvature is strictly positive, we only need to consider the doubly warped product end. For $t\geq R_1$, one can apply Lemma \ref{lmmriccicurvature} to show that
\begin{eqnarray}
\Ric_{\hat{g}} \left(X_0 \right) & = & \frac{\alpha(1-\alpha)}{(t+C)^2} X_0 ,\\
\Ric_{\hat{g}} \left(X_1 \right) & = & \left[ (k-2)\frac{1-(1-\epsilon)^2 }{(1-\epsilon)^2(t+C)^2} - \frac{2\alpha}{(t+C)^2 } \right] X_1 ,\\
\Ric_{\hat{g}} \left(X_2 \right) & = & \left[ \frac{\alpha(1-\alpha)}{(t+C)^2} + \frac{1-\delta^2\alpha^2(t+C)^{2\alpha-2} }{\delta^2(t+C)^{2\alpha}} - \frac{(k-1)\alpha}{(t+C)^2 } \right] X_2 ,
\end{eqnarray}
where $X_0=\frac{\partial}{\partial t}$, $X_1\in T\sphere^{k-1}$ and $X_2\in T\sphere^2$ are unit vector fields. Hence for $\epsilon\in (0,\frac{1}{10})$, $\alpha\in(0, e^{-100}\epsilon)$ and $\delta\in(0,\alpha^{10})$, we have
$$ \Ric_{\hat{g}}\geq \frac{\alpha^{10}}{(t+C)^2}\geq \frac{\alpha^{10}}{2(1+C)^2}\cdot\frac{1}{t^2+1} .$$
This gives the required Ricci lower bound, and completes
the proof.
\end{proof}

\section{Construction of the sequence of metrics}
\label{sectioninduction}
In this section, we construct the sequence of $(k+2)$-dimensional manifolds $(M_i ,g_i )  \stackrel{GH}{\longrightarrow} (X_\epsilon ,d_\epsilon ) $, where $(X_\epsilon ,d_\epsilon )$ is a metric space close to $(X,h)$ in the Gromov-Hausdorff sense, and $(X , h)$ is a $k$-dimensional smooth complete Riemannian manifold with $k\geq 6$ and $\Ric_h > \lambda \in\mathbb{R} $. 

\subsection{The surgery}
Before constructing the sequence of metrics, we need to prove the following lemma. This is the surgery for constructing singularities.

\begin{lmm}
\label{lmm_surgery}
Fix $k\geq 6$. Let $(B,g_B)$ be a $k$-dimensional Riemannian metric (not necessary complete) and let $f_B$ be a strictly positive function on $B$. Let $g$ denote the warped product metric $g=g_B+f_B^2g_{\sphere^2}$ on $B\times \sphere^2$ such that $\Ric_{g}\geq\lambda+\kappa$ for some $\kappa\in(0,1)$. Then we have:
\begin{enumerate}[\rm (i)${}_4$.]
\item A family of Riemannian metrics $g_{B;j}$ on $B$, where $j\in\mathbb{Z}_{\geq 0}$,
\item A family of strictly positive functions $f_{B;j}$ on $B$,
\item A precompact open subset $W\subset \bar{W}\subset B$ diffeomorphic to a ball,
\item A smooth embeddings $\iota_{B}:\bar{\ball}^3_1\to U$, where $\bar{\ball}^3_1=\{y\in\mathbb{R}^3:|y|\leq 1\}$,
\item A sequence of positive constants $\delta_j\in (0,1)$,
\end{enumerate}
such that:
\begin{enumerate}[\rm (i)${}_4$.]
\addtocounter{enumi}{5}
\item The Ricci curvature $\Ric_{g_{j}}>\lambda$, where $g_{j}=g_{B;j}+f_{B;j}^2 g_{\sphere^2}$,
\item The functions $f_{B;j}=\delta_j f_B$ on $B\setminus W$,
\item The metrics $g_{B;j}$ converge to a Riemannian metric on $B\setminus \iota_B(\sphere_1^2)$ smoothly,
\item The restriction of $g_{B;j}|_{\iota_B(\sphere_1^2)}$ converge to $0$ smoothly,
\item The metrics satisfy $g_{B;j}\leq 2g_{B;0}$.
\end{enumerate}
\end{lmm}

\begin{proof}
Let $p\in B$, and let $\epsilon_p>0$ be a small constant such that the ball $W=B_{3\epsilon_p}(p)\subset B$ is precompact, and the rescaled ball $(W,\epsilon_p^{-2}g_B)$ satisfies the condition (\ref{huppnaberwangcondition1})-(\ref{huppnaberwangcondition3}) in Lemma \ref{huppnaberwanglemma}.

Then for any choices of parameters $0 < \eps < \eps_0(n,|\lambda|)$, $0<\alpha < \alpha_0(n,\eps)$, $0 < \hat{r} < \hat r_0(n,\alpha, \lambda,\eps)$, and $0 < \hat{\delta} < \hat{\delta}_0(n,\lambda,\alpha,\epsilon,||f_B||_{L^\infty},\hat{r})$, there exists $0 < \delta = \Psi(\,\hat{\delta}\|f_B\|_\infty\mid n, \alpha, \lambda, \eps)$ and a warped product metric $\hat g = \hat{g}_B+\hat{f}_{B}^{2} g_{\sphere^2}$ such that:
\begin{enumerate}[\rm (a$'$).]
	\item The Ricci lower bound $\Ric_{\hat g}>\lambda+\kappa - C_1(k)\eps_p^{-2}\eps$ holds for $\hat r/2 \leq r \leq 2\eps_p$\, ,
	\item $\hat{g} = g_B+\hat{\delta}^2 f_B^2 g_{\sphere^2}$, $ \forall r\in [1,2] $ ,
	\item $\hat{g} = dr^2+(1-\eps)^2 r^2 g_{\sphere^k}+\delta^2 r^{2\alpha} g_{\sphere^2}$, $\forall r\leq \hat r$ ,
	\item The identity map $\Id:(\ball_2(p),\hat{g}_B)\to (\ball_2(p),g_B)$ is $(1+2\epsilon)$-bi-Lipschitz.
\end{enumerate}

Assume that $\epsilon\leq\frac{\kappa\epsilon_p^2}{10(1+C_1(k))}$. Then $\Ric_{\hat g}>\lambda$.

Now we apply Proposition \ref{prop_eightdim_localmodel} to $p\in W$. After scaling, for any $\varsigma\in(0,1)$, one can obtain a sequence of complete Riemannian metrics $\hat{g}_{\varsigma;j}$ on $\mathbb{R}^k\times\sphere^2$ such that $\Ric_{\hat{g}_{\varsigma;j}}\geq \frac{\tau}{\varsigma^2+r_j^2}$ and there exists a compact subset $U\subset\mathbb{R}^k$ satisfies 
$$\big((\mathbb{R}\setminus U)\times\sphere^2,\hat{g}_{\varsigma;j}\big)\cong\big([\varsigma,\infty)\times\sphere^{k-1}\times\sphere^2,dr_j^2+(1-\eps)r_j^2g_{\sphere^k}+\varsigma^{2-2\alpha}\delta_p^2r_j^{2\alpha}g_{\sphere^2}\big),$$
where $\tau,\delta_p>0$ are constants independent of $j$.

Let $\hat{r},\varsigma>0$ be constants such that $\varsigma^{1-\alpha}\delta_p<\hat{\delta}$ and $\frac{\tau}{\varsigma^2+r_j^2}>\lambda+1$, $\forall r_j\in(0,\hat{r})$. By choosing $\delta=\varsigma^{1-\alpha}\delta_p$, we can glue the metrics constructed in Proposition \ref{prop_eightdim_localmodel} to $B_{\hat{r}}(p)$, and got the sequence of metrics we need. The required properties now follows from Proposition \ref{prop_eightdim_localmodel}.
\end{proof}

\subsection{Induction progress}
Let $X_{0,0}$ denote the Riemannian manifold $(X,h)$. Let $(M_{0,0} ,g_{0,0})$ denotes the Riemannian manifold $\left( X\times \sphere^2 , h+ \frac{e^{-2i}}{1+|\lambda|^2} g_{\sphere^2} \right) $. It follows from $\Ric_h >\lambda $ that $\Ric_{g_{0,0}} > \lambda $.

Our next step is to construct Riemannian manifolds $(M_{N,i} ,g_{N,i} )=\left( X\times \sphere^2 , h_{N,i}+ f_{N,i}^2 g_{\sphere^2} \right)$ and metric spaces $(X_N,d_N)$ by induction.

\begin{lmm}
\label{lmm_induction}
Fix $p_0 \in X$ and $\delta\in(0,1)$. Then we have:
\begin{enumerate}[\rm (I).]
\item A family of complete Riemannian metrics $h_{N,i}$ on $X$, where $N,i\in\mathbb{Z}_{\geq 0}$,
\item A family of strictly positive functions $f_{N,i}$ on $X$,
\item A sequence of smooth maps $\mathscr{P}_{N}:X\to X$,
\item A family of smooth embeddings $\iota_{N}:S_N\times\bar{\ball}^3_1\to X$, where $S_N=\{1,\cdots,N\}\subset\mathbb{N}$, and $\bar{\ball}^3_1=\{y\in\mathbb{R}^3:|y|\leq 1\}$,
\item A family of precompact open neighborhoods $W_{N}$ of $\iota_{N}(\{N\}\times\bar{\ball}^3_1)\subset X$,
\item Three increasing sequences of integers $\{N_\ell\}_{\ell=1}^\infty$, $\{i_\ell\}_{\ell=1}^\infty$, and $\{C_N\}_{N=1}^\infty$
\end{enumerate}
such that:
\begin{enumerate}[\rm (I).]
\addtocounter{enumi}{6}
\item The Ricci curvature of $g_{N,i}=h_{N,i}+f_{N,i}^2g_{\sphere^2}$ satisfies $\Ric_{g_{N,i}}>\lambda$,
\item $\iota_{N+1}|_{S_N\times\bar{\ball}_1^3}=\iota_{N}$, and $\iota_{N}(S_{N-1}\times\bar{\ball}_1^3)\subset X\setminus (W_N\cup\{p_0\})$,
\item The diameter $\mathrm{diam}_{h_{N,i}} (W_{N})+\mathrm{diam}_{h_{N-1,i}} (W_{N})\leq \frac{i+1}{i+2}\cdot2^{-e^{10(C_N+N+1)}}\delta$, $\forall i\in\mathbb{Z}_{\geq 0}$,
\item $\mathscr{P}_{N}(p_0)=p_0$, and $\mathscr{P}_{N}^* h_{N-1,i}\leq C_Nh_{N,i}$,
\item The restriction $\mathscr{P}_{N}|_{X\setminus W_{N}}$ gives the identity map $\Id_{X\setminus W_N}$,
\item The restriction $h_{N-1,i}|_{X\setminus W_{N}}=h_{N,i}|_{X\setminus W_{N}}$, $i\in\mathbb{Z}_{\geq 0}$,
\item The function $|f_{N,i}|\leq\frac{1}{N+1}$ on $X$, and the quotient $f_{N,i}f_{N-1,i}^{-1}$ is a constant on $X\setminus W_N$, $i\in\mathbb{Z}_{\geq 0}$,
\item The pointed Riemannian manifolds $(X,h_{N,i},p_0)$ converge to $(X_N,d_N,p_{N;\infty})$ as $i\to\infty$ in the Gromov-Hausdorff sense,
\item The embeddings $\iota_{N}:S_N\times\bar{\ball}^3_1 \to X$ converge to $\iota_{N,\infty}:S_N\times\bar{\ball}^3_1 \to X_N$ as $i\to\infty$ in the Gromov-Hausdorff sense,
\item The metric space $X_N$ is homeomorphic to quotient space $X/\sim_N$, where $x_1\sim_N x_2$ if and only if $x_1,x_2\in \iota_{N,i}(\{j\}\times\sphere^2_1)$ for some $j\in S_N$,
\item The metrics $h_{N,i}$ and functions $f_{N,i}$ converge smoothly to a Riemannian metric and $h_{N,\infty}$ a smooth positive function $f_{N,\infty}$ on $X_N\setminus \iota_{N,\infty}(S_N\times\sphere_1^2)$ as $i\to\infty$, respectively,
\item The map $\iota_{N,\infty}$ factors through the quotient map $\bar{\ball}^3_1 \to\bar{\ball}^3_1 /\sphere^2_1$,
\item For any $\ell\in\mathbb{N}$, $N\geq N_\ell$ and $p\in B_{e^\ell}(p_0)\subset (X,h_{N,i})$, 
$$\#\left\{j\in S_N: \bar{W}_{j}\subset B_{e^{-\ell}}(p)\right\}\geq \ell,$$
\item For any $\ell\in\mathbb{N}$, $N\leq N_\ell$ and $i\geq i_\ell$, $f_{N,i+1}=f_{N,i}$.
\end{enumerate}
\end{lmm}

\begin{proof}
At first, we describe the construction by induction on $\ell\in\mathbb{Z}_{\geq 0}$. For $\ell=0$, let $N_0=0$, $h_{0,i}=h$, $\mathscr{P}_{0,i}=\Id_X$, $W_0=\varnothing$, $i_0=0$, and $f_{0}=(10+10|\lambda|)^{-10}$. 

Let $\ell\in\mathbb{N}$. Suppose that all the data with $N\leq N_{\ell-1}$ and $i\in\mathbb{Z}_{\geq 0}$ have been constructed. Now we find an integer $N_\ell>N_{\ell-1}$ and construct the next data for $N_{\ell-1}<N\leq N_\ell$. By the induction hypothesis, we can apply the property (XVII) to show that there exists a finite set
$$\mathbf{P}_{\ell}=\{p_{\ell;j}\}_{j=1}^{N'_\ell}\subset \bigcup_{i=0}^\infty \left(B^{h_{N_{\ell-1,i}}}_{e^{10\ell}}(p_0)\right)\bigg\backslash \left(\bigcap_{i=0}^\infty\left(B^{h_{N_{\ell-1,i}}}_{e^{-10\ell}}(p_0)\right)\bigcup\iota_{{N_{\ell-1}}}(S_{N_{\ell-1}}\times \bar{\ball}^3)\right),$$ 
such that for any $i\in\mathbb{Z}_{\geq 0}$ and $p\in B^{h_{N_{\ell-1,i}}}_{e^{5\ell}}(p_0)\setminus B^{h_{N_{\ell-1,i}}}_{e^{-5\ell}}(p_0)$, 
$$\#(\mathbf{P}_{\ell}\cap B^{h_{N_{\ell-1,i}}}_{e^{-10\ell}}(p_0))\geq 100\ell,$$
where $B^{h_{N_{\ell-1,i}}}_r(p)$ is the ball in $(X,h_{N_{\ell-1,i}})$ of radius $r$ centered at $p$.

Similarly, for any $\ell\in\mathbb{N}$, we can find a constant $r_{\ell}\in(0,e^{-100e^{100(N_{\ell-1}+N'_\ell)}})$ such that for any $i\in\mathbb{Z}_{\geq 0}$, the balls $\left\{B^{h_{N_{\ell-1,i}}}_{2^{100}r_{\ell}}(p_{\ell;j})\right\}_{j=1}^{N'_\ell}$ are pairwise disjoint, $B^{h_{N_{\ell-1,i}}}_{2^{100}r_{\ell}}(p_{\ell;j})\cap\iota_{{N_{\ell-1}}}(S_{N_{\ell-1}}\times \bar{\ball}^3)=\varnothing$, and the renormalized domains $\left(B^{h_{N_{\ell-1,i}}}_{2^{100}r_{\ell}}(p_{\ell;j})\times\sphere^2,2^{-200}r_\ell^{-2}h_{N_{\ell-1,i}}+f_{N_{\ell-1},i}^2g_{\sphere^2}\right)$ satisfying the conditions in Lemma \ref{huppnaberwanglemma}. Set 
$$N_\ell=N_{\ell-1}+N'_\ell,\quad \textrm{ and }\quad U_{N}=\left(\bigcap_{i=1}^\infty B^{h_{N_{\ell-1,i}}}_{100\delta r_{\ell}}(p_{\ell;N-N_{\ell-1}})\right)^{\mathrm{int}},$$
where $Z^{\mathrm{int}}$ is the interior of $Z$. It is easy to see that $U_N\neq\varnothing$ is an open neighborhood of $p_{\ell,N-N_{\ell-1}}$.

Now we construct $\mathscr{P}_{N}$ and $g_{N,i}$ for $N_{\ell-1}+1\leq N\leq N_\ell$ by induction on $N$. Suppose that all the data with $N\leq \hat{N}\in\mathbb{Z}\cap[N_{\ell-1}, N_\ell-1]$ and $i\in\mathbb{Z}_{\geq 0}$ have been constructed.

By the property (XVII) again, we can assume that there are coordinates $\mathfrak{F}_N=(u_{N,1},\cdots,u_{N,k}):U_N\to\mathbb{R}^k$ such that $\mathfrak{F}_N(p_{\ell,N-N_{\ell-1}})=0$, the image $\mathfrak{F}_N(U_N)$ contains an Euclidean ball $B_{2^{100}\delta r'_\ell} (0)\subset\mathbb{R}^k$, and for any $i$, $e^{-c_N}h_{\Euc,N}\leq h_{N_{\ell-1},i}\leq e^{c_N}h_{\Euc,N}$, where $r'_\ell,c_N>0$ are constants, and $h_{\Euc,N}=\sum_{j=1}^k du_{N,j}\otimes du_{N,j}$ is the Euclidean metric on $\mathbb{R}^k$. Let $\mathscr{P}_{N;\Euc}:B_{2^{100}\delta r'_\ell}(0)\to B_{2^{100}\delta r'_\ell}(0)\subset\mathbb{R}^k$ be a smooth map such that $\mathscr{P}_{N;\Euc}(B_{r'_\ell}(0))=0$, $\mathscr{P}_{N;\Euc}|_{B_{2^{100}\delta r'_\ell}(0)\setminus B_{2^{90}\delta r'_\ell}(0)}=\Id_{B_{2^{100}\delta r'_\ell}(0)\setminus B_{2^{90}\delta r'_\ell}(0)}$, and $\|D\mathscr{P}_{N;\Euc}\|_{h_{\Euc,N}}\leq 2^{10k}$. Let $\mathscr{P}_{N}=\mathfrak{F}_N^{-1}\circ\mathscr{P}_{N;\Euc}\circ\mathfrak{F}_N$. 

Let $W_N\subset U_N$ denote an open neighborhood of $p_{\ell,N-N_{\ell-1}}$ defined by $\mathfrak{F}_N^{-1}(B_{2^{95}\delta r'_\ell}(0))\subset U_N$, and let $V_N=\mathfrak{F}_N^{-1}(B_{2^{-e^{e^{1000C_N}}}\delta r'_\ell}(0))\subset U_N$, where $C_N\in\mathbb{N}$ is an integer such that $C_N\geq e^{10(c_N+2^{100k}+C_{N-1}+1)}$. For any $i\in\mathbb{Z}$, we can construct metrics $g_{N,i}=h_{N,i}+f^2_{N,i}g_{\sphere^2}$ satisfying (VII), (IX), (XII), (XIII) and (XVI) by applying Lemma \ref{lmm_surgery} to $V_N\times \sphere^2$. Let $\iota_{V_N}$ be the smooth embedding $\iota_{V_N}:\bar{\ball}_1^3\to V_N$ in Lemma \ref{lmm_surgery}. Now we can define $\iota_{N}$ by $\iota_{N}|_{S_{N-1}\times\bar{\ball}_1^3}=\iota_{N-1}$ and $\iota_N|_{\{N\}\times\bar{\ball}_1^3}=\iota_{V_N}$. Then by Lemma \ref{lmm_surgery} and the hypothesis, $\iota_N$ satisfies (VIII), (XV), (XVII), (XVIII) and (XIX). By the argument above, $\mathscr{P}_{N}$ satisfies (X) and (XI).

By Lemma \ref{lmm_surgery} and the property (XVII), for each $N\in\mathbb{N}$, there exists $i_N>0$ such that if $i,i'\geq i_N$, then (VII) holds when we replace $f_{N,i}$ by $f_{N,i'}$. Without loss of generality, we can assume that $f_{N,i}=f_{N,i_N}$, $\forall i\geq i_N$. Hence the proof is completed by letting $i_\ell=\sum_{N\leq N_\ell}i_N$.
\end{proof}

As a corollary of Lemma \ref{lmm_induction}, we have the following estimates.

\begin{coro}
\label{coro_estimate_induction}
Fix $p_0\in X$ and $\delta\in(0,1)$. Let $(X\times \sphere^2,g_{N,i}=h_{N,i}+f^2_{N,i}g_{\sphere^2})$ and $\mathscr{P}_N$ be the data constructed in Lemma \ref{lmm_induction}. Then for any $N,m\in\mathbb{N}$ and  $i\in\mathbb{Z}_{\geq 0}\cup \{\infty\}$, the map 
$$\mathscr{P}^{N+m}_N=\mathscr{P}_{N+m-1}\circ\cdots\mathscr{P}_{N+1}\circ\mathscr{P}_{N}:(X,h_{N+m,i})\to (X,h_{N,i}) $$ 
is a Gromov-Hausdorff $2^{-N-1}\delta$-approximation. Moreover, for any $x_1,x_2\in X$,
$$\dist_{h_{N+m,i}}\left(\mathscr{P}^{N+m}_N(x_1),\mathscr{P}^{N+m}_N(x_2)\right)\leq 2\cdot\max\left\{\dist_{h_{N,i}} (x_1,x_2),\dist_{h_{N,i}} (x_1,x_2)^{1-\frac{N+m}{2(N+m+1)}}\right\} .$$
\end{coro}

\begin{proof}
Let $W_N$ be the subsets constructed in Lemma \ref{lmm_induction}. Then we have $\mathscr{P}_N|_{X\setminus W_N}=\Id_{X\setminus W_N}$ and $\mathrm{diam}_{h_{N,i}} (W_{N})+\mathrm{diam}_{h_{N-1,i}} (W_{N})\leq \frac{i+1}{i+2}\cdot2^{-e^{10(N+1)}}\delta$. It follows that for any $N\in\mathbb{N}$ and $a\in\mathbb{Z}_\geq 0$, the map $\mathscr{P}_{N+a}:(X,h_{N+a+1,i})\to (X,h_{N+a,i}) $ is a Gromov-Hausdorff $2^{-10(N+a+1)}\delta$-approximation. Hence $\mathscr{P}^{N+m}_N $ is a Gromov-Hausdorff $2^{-N-1}\delta$-approximation.

Let $C_N>0$ be the constant given in the property (VI) in Lemma \ref{lmm_induction}. When $\dist_{h_{N,i}} (x_1,x_2)\leq 2^{-e^{5(N+C_N+1)}} $, we have 
\begin{eqnarray*}
\dist_{h_{N+1,i}}\left(\mathscr{P}_N(x_1),\mathscr{P}_N(x_2)\right)&\leq&  C_N\cdot \dist_{h_{N,i}} (x_1,x_2) \\
&\leq& C_N\cdot e^{-2^{5(N+C_N+1)}\cdot 2^{-N-5}}\cdot \dist_{h_{N,i}} (x_1,x_2)^{1-2^{-N-5}}\\
&\leq & \dist_{h_{N,i}} (x_1,x_2)^{1-2^{-N-5}}.
\end{eqnarray*}

If $\dist_{h_{N,i}} (x_1,x_2)\geq 2^{-e^{5(N+C_N+1)}} $, we can conclude that
\begin{eqnarray*}
\dist_{h_{N+1,i}}\left(\mathscr{P}_N(x_1),\mathscr{P}_N(x_2)\right)&\leq&  \dist_{h_{N,i}} (x_1,x_2)+2^{10-e^{10(C_N+N+1)}}\delta \\
&\leq& \left(1+ 2^{10-e^{10(C_N+N+1)}} \cdot 2^{e^{5(N+C_N+1)}}\cdot \delta\right)\cdot \dist_{h_{N,i}} (x_1,x_2) \\
&\leq & \left(1+ 2^{-e^{8(C_N+N+1)}} \right)\cdot \dist_{h_{N,i}} (x_1,x_2) .
\end{eqnarray*}

It follows that
\begin{eqnarray*}
    & & \dist_{h_{N+1,i}}\left(\mathscr{P}_N(x_1),\mathscr{P}_N(x_2)\right)\\
    &\leq& \left(1+ 2^{-e^{8(N+1)}} \right)\cdot \max\left\{\dist_{h_{N,i}} (x_1,x_2),\dist_{h_{N,i}} (x_1,x_2)^{1-2^{-N-5}} \right\},
\end{eqnarray*}
and hence
\begin{eqnarray*}
    & & \dist_{h_{N+m,i}}\left(\mathscr{P}^{N+m}_N(x_1),\mathscr{P}^{N+m}_N(x_2)\right)\\
    &\leq& \prod_{a=0}^{m-1}\left(1+ 2^{-e^{8(N+a+1)}} \right)\cdot \max\left\{\dist_{h_{N,i}} (x_1,x_2),\dist_{h_{N,i}} (x_1,x_2)^{\prod_{a=0}^{m-1}(1-2^{-N-a-5})} \right\}\\
    &\leq& 2\cdot\max\left\{\dist_{h_{N,i}} (x_1,x_2),\dist_{h_{N,i}} (x_1,x_2)^{1-\frac{N+m}{2(N+m+1)}}\right\} ,
\end{eqnarray*}
and the corollary is proved.
\end{proof}

\section{Proof of Theorem \ref{thmnonmfd}}
\label{sectionproof}
In this section, we will prove Theorem \ref{thmnonmfd}. For convenience, we recall the statement of Theorem \ref{thmnonmfd} here.

\begin{thm}[Imply Theorem \ref{thmnonmfd}]
\label{thmnonmfd2}
Let $(X , h)$ be an $k$-dimensional smooth complete Riemannian manifold with $k\geq 6$ and $\Ric_h > \lambda \in\mathbb{R} $. Then for each $\epsilon >0$, there exists a metric space $(X_\epsilon ,d_\epsilon )$ such that
\begin{enumerate}[\rm (i).]
\item $d_{GH} (X,X_\epsilon ) \leq \epsilon $,  \\[-3mm]
\item $(X_\epsilon ,d_\epsilon ) $ is $k$-rectifiable, \\[-3mm]
\item There exists a sequence of $(k+2)$-dimensional boundary free complete manifolds $(X\times\mathbb{S}^2 ,g_i ) \stackrel{GH}{\longrightarrow} (X_\epsilon ,d_\epsilon )$ with $\Ric_{g_i} > \lambda $,  \\[-3mm]
\item For any nonempty open subset $U$ of $X_\epsilon $, there exists a continuous map $\mathfrak{P}:\mathfrak{N}_\ell\times \sphere^3\to U$ such that the homomorphism $\mathfrak{P}_*:H_3 (\mathfrak{N}_\ell\times \sphere^3;\mathbb{Z}) \to H_3 (X_\epsilon ;\mathbb{Z})$ induced by $\iota$ is an injection, where $\mathfrak{N}_\ell=\{1,\cdots,\ell\}$.
\end{enumerate}
\end{thm}

\begin{rmk}
Since $H_3 (\mathfrak{N}_\ell\times \sphere^3;\mathbb{Z})\cong\mathbb{Z}^\ell$, the property (iv) of Theorem \ref{thmnonmfd2} implies that the local $3$-homology group ${\rm{im}} \left( H_3 (U;\mathbb{Z}) \to H_3 (X_\epsilon ;\mathbb{Z}) \right) $ are infinitely generated. Consequently, every open subset of $X_\epsilon $ is not a topological manifold.
\end{rmk}

\begin{proof}
Let $\delta=2^{-100}\epsilon$. By Lemma \ref{lmm_induction}, there exists a family of Riemannian manifolds $(M_{N,i},g_{N,i})=(X\times \sphere^2,h_{N,i}+f_{N,i}^2g_{\sphere^2})$ such that 
$$d_{GH} ((X,h),(X_{N,i},h_{N,i}) )\leq 10\cdot(1-2^{-N-1})\cdot\delta \leq \epsilon.$$
It follows that for any $N\in\mathbb{N}$, $d_{GH} ((X,h),(X_N,d_N) )\leq\epsilon$. Since $\{(X_N,d_N)\}_{N\in\mathbb{N}}$ is a Cauchy sequence in the Gromov-Hausdorff sence, we can conclude that this sequence converges to a metric space $(X_\infty,d_\infty)$ in the Gromov-Hausdorff sence, and $d_{GH} (X,X_\infty ) \leq \epsilon $. Moreover, by an argument similar to \cite{hnw1} and \cite{zhou1}, $X_\infty $ is $k$-rectifiable. By the diagonal argument, we can find a sequence $\{M_{N_j,i_j}\}_{j=1}^\infty$ converging to $(X_\infty,d_\infty)$. Note that $|f_{N,i}|\leq \frac{1}{N+1}$ implies that $d_{GH}\big((X_{N,i},h_{N,i}),(M_{N,i},g_{N,i})\big)\leq\frac{10}{N+1}$. By choosing $X_\epsilon=(X_\infty,d_\infty)$, we obtain (i)-(iii).

Now we only need to prove (iv). Let $U\subset X_\infty$ be an open subset. Then there exists a sequence of points $q_i\in (X_N,d_N)$ converging to $q_\infty \in U\subset X_\infty$. Let $r_\infty>0$ be a constant such that $B_{10r_\infty}(p_\infty)\subset U$. 

By Corollary \ref{coro_estimate_induction}, for any $N\in\mathbb{N}$, the maps $\mathscr{P}_N^{N+m}$ converging to a H\"older continuous map $\mathscr{P}_N^{\infty}:X_\infty\to X_N$ as $m\to\infty$ in the Gromov-Hausdorff sense, and $\mathscr{P}_N^{\infty}$ is a Gromov-Hausdorff $2^{-N}\delta$-approximation. By Lemma \ref{lmm_induction}, there exists a map $\iota_{\infty,\infty}:\mathbb{N}\times\sphere^3\to X_\infty$ such that for any $N\in\mathbb{N}$, the restriction $\iota_{\infty,\infty}|_{S_N\times\sphere^3}=\iota_{N,\infty}:S_N\times\sphere^3\to X_N$, where $S_N=\{1,\cdots,N\}$, and $\iota_{N,\infty}$ is the map in Lemma \ref{lmm_induction} given by the limit of the sequence of smooth embeddings $\iota_{N,i}$. Apply Lemma \ref{lmm_induction} again, one can see that there exists an integer $N_{q_\infty;\ell}\in\mathbb{N}$ such that for any $N\geq N_{q_\infty;\ell}$,
$$\#\left\{j\in S_N: \iota_{N,\infty}(\{j\}\times\sphere^3)\subset B_{r_\infty}(q_N) \right\} \geq 2\ell .$$

Fix an integer $N\geq N_{q_\infty;\ell}$. Now we embed $\mathfrak{N}_\ell $ into the set $\left\{j\in S_N: \iota_{N,\infty}(\{j\}\times\sphere^3)\subset B_{r_\infty}(q_N) \right\}$, and define $\mathfrak{P}$ by the restriction $\mathfrak{P}=\iota_{\infty,\infty}|_{\mathfrak{N}_\ell}$. By considering the composition $\mathscr{P}^\infty_{N}\circ \mathfrak{P}$, we only need to show that for any $N\in\mathbb{N}$, the homeomorphism
$$(\iota_{N,\infty})_*:H_3 (S_N\times \sphere^3;\mathbb{Z}) \to H_3 (X_N ;\mathbb{Z}) $$
induced by $\iota_{N,\infty}:S_N\times\sphere^3\to X_N$ is an injection.

By construction, $X_N$ is homotopic equivalent to the glued space $X\cup_{\iota_N(S_N\times \sphere^2_1)}(S_N\times \bar{\ball}^3_1)$. Note that $\bar{\ball}^3_1$ can be seen as the topological cone over $\sphere^2$. Since $\bar{\ball}^3_1$ is contractible, $X_N$ is homotopic equivalent to $X\vee\bigvee_{j\in S_N} \sphere^3$, the wedge sum. Hence $H_3 (X_N ;\mathbb{Z})\cong H_3\big(X\vee\bigvee_{j\in S_N} (\{j\}\times\sphere^3);\mathbb{Z}\big)$. Moreover, the homological class $(\iota_{N,\infty})_*( e_{\{j\}\times\sphere^3})\in H_3 (X_N ;\mathbb{Z})$ is corresponding to $\mathrm{pr}_*(e_{\{j\}\times\sphere^3})$, where $e_{\{j\}\times\sphere^3}\in H_3(\{j\}\times\sphere^3)$ is the generator, and $\mathrm{pr}:S_N\times\sphere^3\to X\vee\bigvee_{j\in S_N} (\{j\}\times\sphere^3)$ is the canonical map. By the Mayer–Vietoris sequence, we can conclude that $\mathrm{pr}_*:H_3(S_N\times\sphere^3)\to H_3\big(X\vee\bigvee_{j\in S_N} (\{j\}\times\sphere^3)\big)$ is an injection. It follows that $(\iota_{N,\infty})_*$ is also an injection, and the proof is completed.
\end{proof}


\begin{thebibliography}{BBEGZ16} 
{

\bibitem{bd1} E. Bru\`e, Q. Deng:
\emph{Topology of two-dimensional collapsed spaces with lower Ricci bounds},
preprint, arXiv: 2606.19189.

\bibitem{burdick1} B. L. Burdick:
\emph{Metrics of positive Ricci curvature on the connected sums of products with arbitrarily many spheres},
Ann. Global Anal. Geom. $\mathbf{58}$ (2020), 433-476.

\bibitem{chco3} J.~Cheeger, T. H. Colding:
\emph{On the structure of spaces with Ricci curvature bounded below. I},
J. Differential Geom. $\mathbf{46}$ (1997), 406–480.

\bibitem{chco4} J.~Cheeger, T. H. Colding:
\emph{On the structure of spaces with Ricci curvature bounded below. II},
J. Differential Geom. $\mathbf{54}$ (2000), 13–25.

\bibitem{chco5} J.~Cheeger, T. H. Colding:
\emph{On the structure of spaces with Ricci curvature bounded below. III},
J. Differential Geom. $\mathbf{54}$ (2000), 37–74.

\bibitem{jcwsjan1} J.~Cheeger, W. Jiang, A.~Naber:
\emph{Rectifiability of singular sets of noncollapsed limit spaces with Ricci curvature bounded below},
Ann. of Math. $\mathbf{193}$ (2021), 407–538.

\bibitem{chti1} J.~Cheeger, G. Tian:
\emph{Curvature and injectivity radius estimates for Einstein $4$-manifolds},
J. Amer. Math. Soc. $\mathbf{19}$ (2006), 487-525.

\bibitem{colding1} T. H. Colding:
\emph{Ricci Curvature and Volume Convergence},
Ann. of Math. $\mathbf{145}$ (1997), 477-501.

\bibitem{coldingnaber1} T. H. Colding, A.~Naber:
\emph{Sharp H\"older continuity of tangent cones for spaces with a lower Ricci curvature bound and applications},
Ann. of Math. $\mathbf{176}$ (2012), 1173–1229.

\bibitem{lina1} L. Chen: 
\emph{A remark on regular points of Ricci limit spaces},
Front. Math. China $\mathbf{11}$ (2016), 21-26.

\bibitem{dephi1} G. De Philippis, N. Gigli:
\emph{Non-collapsed spaces with Ricci curvature bounded from below},
J. \'Ec. Polytech. Math. $\mathbf{5}$ (2018), 613–650.

\bibitem{gigli1} G. De Philippis, N. Gigli:
\emph{On the differential structure of metric measure spaces and applications},
Mem. Amer. Math. Soc. $\mathbf{236}$ (2015), no. 1113, vi+91 pp.

\bibitem{geji1} H. Ge, W. Jiang:
\emph{$\epsilon$-regularity for shrinking Ricci solitons and Ricci flows},
Geom. Funct. Anal. $\mathbf{27}$ (2017), 1231-1256.

\bibitem{aha1} A.~Hatcher:
\emph{Algebraic topology},
Cambridge University Press, Cambridge, 2002.

\bibitem{honda1} S. Honda: 
\emph{On low-dimensional Ricci limit spaces},
Nagoya Math. J. $\mathbf{209}$ (2013), 1-22.

\bibitem{hnw1} E. Hupp, A. Naber, K.-H. Wang: 
\emph{Lower Ricci Curvature and Nonexistence of Manifold Structure},
Geom. Topol. $\mathbf{29}$ (2025), 443–477.

\bibitem{km1} V. Kapovitch, A. Mondino:
\emph{On the topology and the boundary of $N$-dimensional $\mathrm{RCD} (K,N)$ spaces},
Geom. Topol. $\mathbf{25}$ (2021) 445–495.

\bibitem{kw1} V. Kapovitch, B. Wilking:
\emph{Structure of fundamental groups of manifolds with
 Ricci curvature bounded below},
preprint, arxiv: 1105.5955.

\bibitem{mrs1} D. Maximo, P. Reiser, D. Semola:
\emph{Ricci curvature and minimal hypersurfaces with large Betti numbers},
Calc. Var. $\mathbf{65}$ (2026).

\bibitem{menguy2} X. Menguy:
\emph{Noncollapsing examples with positive Ricci curvature and infinite topological type},
Geom. Funct. Anal. $\mathbf{10}$ (2000), 600–627.

\bibitem{naber1} A. Naber:
\emph{Conjectures and open questions on the structure and regularity of spaces with lower ricci curvature bounds},
SIGMA Symmetry Integrability Geom. Methods Appl. $\mathbf{16}$ (2020), 8 pp.

\bibitem{panwei1} J. Pan, G. Wei:
\emph{Examples of Ricci limit spaces with non-integer Hausdorff dimension},
Geom. Funct. Anal. $\mathbf{32}$ (2022), 676–685.

\bibitem{pp1} P.~Petersen:
\emph{Riemannian geometry},
3rd ed., Grad. Texts in Math, 171. Springer, Cham, 2016.

\bibitem{perel1} G.~Perelman:
\emph{Construction of manifolds of positive Ricci curvature with big volume and large Betti numbers},
Comparison geometry (Berkeley, CA, 1993-94), 157-163, Math. Sci. Res. Inst. Publ., 30, Cambridge Univ. Press, Cambridge, 1997.

\bibitem{reiser1} P. Reiser:
\emph{Metrics of Positive Ricci Curvature on Simply-Connected Manifolds of Dimension $6k$},
J. Topol. $\mathbf{17}$ (2024), Paper No. e70007, 50 pp.

\bibitem{jkw1} J. Wang: 
\emph{Ricci Limit Spaces Are Semi-locally Simply Connected},
J. Differential Geom. $\mathbf{128}$ (2024), 1301-1314.

\bibitem{zhou2} S. Zhou: 
\emph{On the Gromov-Hausdorff limits of Tori with Ricci conditions},
preprint, arXiv: 2309.10997. Accepted by Amer. J. Math..

\bibitem{zhou1} S. Zhou: 
\emph{Examples of Ricci limit spaces with infinite holes},
preprint, arXiv: 2404.00619.

}\end{thebibliography}
\end{document}